\documentclass[pdflatex,sn-mathphys]{sn-jnl}% Math and Physical Sciences Reference Style
\usepackage{amsmath}
\usepackage{ulem}
\usepackage{pdfpages}
\usepackage{lineno}
\jyear{2023}%
\usepackage{url}
\usepackage{amsmath}
\theoremstyle{thmstyleone}%
\theoremstyle{thmstyletwo}%

\theoremstyle{thmstylethree}%

\let\oldequation\equation
\let\oldendequation\endequation

\renewenvironment{equation}
  {\linenomathNonumbers\oldequation}
  {\oldendequation\endlinenomath}

\newcommand{\ed}[1]{{\color{black} #1}}
\newcommand{\paul}[1]{{\color{black} #1}}
\newcommand{\anil}[1]{{\color{black} #1}}
\begin{document}

\title[Fractal Basins as a Mechanism for the Nimble Brain]{Fractal Basins as a Mechanism for the Nimble Brain}
%\linenumbers

\author*[1,2]{\fnm{Erik} \sur{Bollt}}\email{bolltem@clarkson.edu}

\author[1,2]{\fnm{Jeremie} \sur{Fish}}\email{jafish@clarkson.edu}

\author[1,2]{\fnm{Anil} \sur{Kumar}}\email{anikumar@clarkson.edu}

\author[1,2,3]{\fnm{Edmilson} \sur{Roque dos Santos}}\email{edmilson.roque.usp@gmail.com}

\author[4]{\fnm{Paul J.} \sur{Laurienti}}\email{plaurien@wakehealth.edu}

\affil*[1]{\orgdiv{Department of Electrical and Computer Engineering}, \orgname{Clarkson University}, \orgaddress{\street{8 Clarkson Ave.}, \city{Potsdam}, \postcode{13699}, \state{New York}, \country{U.S.A}}}

\affil[2]{\orgdiv{Clarkson Center for Complex Systems Science}, \orgname{Clarkson University}, \orgaddress{\street{8 Clarkson Ave.}, \city{Potsdam}, \postcode{13699}, \state{New York}, \country{U.S.A}}}

\affil[3]{\orgdiv{Instituto de Ci\^encias Matem\'aticas e Computa\c{c}\~ao}, \orgname{Universidade de S\~ao Paulo}, \orgaddress{\street{Av. Trab. São Carlense, 400 }, \city{S\~ao Carlos}, \postcode{ 13566-590}, \state{S\~ao Paulo}, \country{Brazil}}}

\affil[4]{\orgdiv{Department of Radiology}, \orgname{Wake Forest University School of Medicine}, \orgaddress{\street{475 Vine Street}, \city{Winston-Salem}, \postcode{27101}, \state{North Carolina}, \country{U.S.A}}}


\abstract{An interesting feature of the brain is its ability to respond to disparate sensory signals from the environment in unique ways depending on the environmental context or current brain state. In dynamical systems, this is an example of multi-stability, the ability to switch between multiple stable states corresponding to specific patterns of brain activity/connectivity. In this article, we describe chimera states, which are patterns consisting of mixed synchrony and incoherence, in a \paul{brain-inspired dynamical systems model} \ed{composed of a network with weak individual interactions and chaotic/periodic local dynamics. We illustrate the mechanism using synthetic time series \paul{interacting on a realistic anatomical brain} network derived from \paul{human} diffusion tensor imaging (DTI)}. We introduce the so-called Vector Pattern State (VPS) as an efficient way of identifying chimera states and mapping basin structures. Clustering similar VPSs for different initial conditions, we show that coexisting attractors of such states reveal intricately ``mingled" fractal basin boundaries that are immediately reachable. This could explain the nimble brain's ability to rapidly switch patterns between coexisting attractors.}


\maketitle

\section{Introduction}\label{sec1}

It is known that the complex dynamics of the brain exhibits numerous spatiotemporal patterns associated with its many capable responses to a given stimulus, as seen in various imaging techniques. Yet, there has not been a good theory to explain how the system is able to switch among these patterns.
Rapidly changing patterns of active brain regions, each containing different types of interconnected neurons that have continuously changing electrochemical properties and environments, only begins to touch on the complexity of a full-scale brain model. This challenge is often countered by course-graining the system to reduce the dimensionality and simplify the model. For instance, instead of analyzing the brain at the neuronal level, even the observational scale of tens of thousands of voxels containing blood oxygenation level dependent (BOLD \cite{ogawa1990brain}) signals from functional magnetic resonance images (fMRI) are down sampled to many fewer anatomical or functional brain regions so that functional brain networks of smaller sizes can be analyzed  \cite{tzourio2002automated, yeo2011organization}.
 
Experiments using fMRI and other imaging technologies reveal that the brain exhibits a rich variety of activity patterns. While it is generally accepted that certain brain regions are more, or less, active when specific tasks are performed or certain sensory systems %\textcolor{brown}{signals} 
such as vision, hearing, or touch are stimulated, %\textcolor{brown}{processed}
it is the global activity patterns that are particularly of interest to us here. An active brain region also implies active neurons, which 
%\sout{, and active} \textcolor{brown}{as} neurons 
share information with other neurons and other brain regions. 
%\sout{Neurons} \textcolor{brown}{They} 
They transmit their signals along axonal pathways via electrical events called action potentials and communicate with other neurons through diverse electrical and chemical synapses \cite{kandel2021principles}. Neural transmission, %\sout{this} \textcolor{brown}{the} 
the process of sharing information along constrained neuroanatomic pathways, can result in neurons exhibiting synchronous large-scale firing patterns, for instance, the collective firing of  neurons generating cortical oscillations \cite{BuzsakiG_2006_complete}. In order to understand how the brain processes environmental cues to generate our experiences, thoughts, and/or emotions it is essential that we better understand these ever-changing, i.e. dynamical patterns of synchronous brain activity \cite{BuzsakiG_2006_complete}.

%%%jjj
%%% k means and a new way here to map synchronous states and coexisting attractor and importantly for brain.
%%% closely packed has the crucial implication here to a mechanisms for...
%%this is the first work that maps fractal basin boundaries in general of complex networked dynamical system and in particular a new way of associating synchronous based chimera states.  ...and in particular for neuro example

 Brain activity can be described mathematically as a complex networked dynamical system which exhibits a key property of multi-stability between numerous states, each associated with different patterns of synchronous activity. The burgeoning field of network neuroscience has used functional brain connectivity \cite{bullmore2009complex} to identify regions of synchronous brain activity, typically assessed using correlations, to show that various patterns of synchrony are associated with distinct cognitive processes \cite{michon2022person, cole2014intrinsic, salehi2020individualized} or brain disorders \cite{wu2021clinical, zhang2021have}. Epilepsy, for example, might be understood as a neurological disease of excess synchrony \cite{Jiruska2013}. Most of the time the brain exhibits patchy or partial synchrony, which is a state in which a subset of nodes (or brain regions) synchronizes while activity in other nodes is incoherent 
 \cite{Schöll_2021}. This state of partial synchrony is often referred to as a chimera state, including cluster synchronization \cite{kaneko1990clustering,belykh2001cluster,abrams2004chimera}. We use the term chimera state broadly to describe the presence of coexisting synchronous and asynchronous (meaning disordered) patterns, and saving ourselves the issue of modifiers to allow for various kinds of synchrony in the definition, see details in the SI. Thus, we consider chimera states as an {\bf activity pattern} where some subset of the system is synchronous and the rest may be incoherent \cite{Majhi2019}. 

\ed{
Chimera states have been observed in brain networks at various scales\paul{, from small to moderate size neural networks composed of spiking neurons \cite{Majhi2019} to brain networks from C. elegans and cats} \cite{hizanidis2016chimera,santos2017chimera}. More recently, researchers have extended their investigations to analyze large-scale functional patterns of \paul{simulated brain activity using various oscillator models interacting on DTI structural brain} networks \cite{Bansal_2019,Kang2019,li2022basins}. Spatiotemporal activity patterns over different brain regions fluctuate over time during resting state, so describing brain dynamics in terms of chimera states holds promise, particularly concerning the \paul{ multistability and metastability of brain activity patterns} \cite{Tognoli2014,Deco2015}.} The key feature of the litany of potential chimera states is that, in a healthy brain, the different organized and disorganized activity patterns coexist with the potential for rapid switching between various states in response to \ed{stimuli}.

\ed{
\noindent
\textbf{Mechanism for the nimble brain.}} \ed{ It has been previously observed that the brain is capable of relatively fast task switching and this has been suggested, with both experimental 
and numerical support \cite{ueltzhoffer2015stochastic,shine2016dynamics,armbruster2012prefrontal,ashourvan2017energy,li2019transitions,loh2007dynamical} to be related to the stability of the basins of attraction involved. Yet, \paul{the dynamical mechanisms} that underpins the ability of the brain to perform such switching in a rapid manner \paul{remain unknown}. In particular, why does the basin of attraction of a particular task appear to be quite stable when it is being performed, while simultaneously allowing for ease of switching between tasks?} In this work, we propose a \ed{potential} mechanism for the agile switching between brain activity patterns/states, a process that supports the nimble brain. Using a perspective of dynamical systems, the nimble brain is explained by a complex basin of attraction for each chimera state with multiple states highly intermingled into a fractal basin boundary. Fractal basin boundaries generally involve a large uncertainty in the final state of a multi-stable system \cite{mcdonald1985fractal}. That is, which initial conditions will lead to a particular final state depends on the detailed intricacies of closely packed and intermingled sets associated with disparate basins of \ed{attraction} \cite{mcdonald1985fractal,grebogi1983fractal,mandelbrot1967long,Mandelbrot_1977_fractals_complete,tel1988fractals}. In particular, there is an apparently rich ``intermingling'' of these boundaries, as the present phenomenon of what is called riddled basins \cite{alexander1992riddled,ott1994transition,  cazelles2001dynamics}, that we present in the results.
This offers a potential mechanism for agile switching between disparate but complex dynamical patterns, i.e. nimble brain activity, because small changes in current state caused by environmental \ed{stimuli} would be enough to switch between distinct stable brain states. 

\ed{An accurate model for capturing the dynamics of the whole-brain has been elusive \cite{Deco2015} \paul{and even if such a model existed, it would be premature to use such a complex, high-dimensional system to map the basin structures investigated here}. Hence, we adopt a simplified model \paul{of spiking neurons on a structural brain network generated using DTI data from a prior study \cite{bonilha2015reproducibility}. Much like prior neuroscience research modeling chimera states \cite{Bansal_2019,Kang2019,li2022basins}, we located brain-inspired dynamical models, Hindmarsh-Rose (HR) neurons in our case, at each node in the DTI network. As a recent research has demonstrated that when coupled, they can exhibit chimera states under specific parameter settings \cite{hizanidis2016chimera}. Others have used models such as Wilson-Cowan oscillators \cite{Bansal_2019,li2022basins}, FitzHugh-Nagumo neurons \cite{Kang2019}, as well as Kuramoto oscilators \cite{li2022basins}. Regardless of the chosen neural model,} this approach allows us to minimize computational complexity while \ed{still} providing a mechanism to emulate the \ed{ essential features of the} nimble brain's behavior. Furthermore, we assess the robustness and general applicability of our findings by testing various individual node dynamics, including Kuramoto oscillators and Hénon maps.} 

We map regions of stability of chimera states to allow us a better understanding of how these disparate patterns co-exist. To make it possible we introduce a technical innovation called the Vector Pattern State (VPS) that characterizes generalized synchronous behaviour from multivariate time series, allowing for phase and approximate synchronization. Using the VPS technology we are able to cluster similar states from different initial conditions and uncover the underlying riddled basin structure of \paul{our brain model}. This observation \ed{sheds light on} a biologically important assertion: fine-scale topological structure of the basins of coexisting chimera states is \paul{potentially underlies the ability of our nimble brain to rapidly switch between various spatial synchronization patterns.} 

\section{Results}\label{sec2}

\subsection{Neuronal model and brain regions}

\ed{Our phenomenological approach is to leverage the presence of chimera states in neuronal systems as a simplified, yet neurologically relevant, model to illustrate our claims regarding the topological fractal basin boundaries in the brain model dynamics. First,} we illustrate the concept of how the brain could switch between disparate pattern states with a semi-synthetic complex coupled system consisting of the well-accepted \paul{HR} model of \ed{spiking neurons}, where the coupling structure is a true \paul{structural brain network with 83 cortical regions connected by white matter fiber tracts measured using DTI.} Fig. \ref{StoryBoard} \paul{illustrates the the organization of this network in brain space}.

A general model of coupled identical units is given by:
\begin{equation}\label{generalcouple}
\mathbf{\dot{x}}_i = f(\mathbf{x}_i) + \sigma \sum \limits_{j = 1}^N [A]_{i,j}h(\mathbf{x}_i,\mathbf{x}_j),
\end{equation}
where $\mathbf{x}_i \in \mathbb{R}^d$ is the state vector, $f: \mathbb{R}^d \to \mathbb{R}^d$ represents the individual node dynamics, $\sigma \in \mathbb{R}^+$ is the coupling strength, $A$ is the adjacency matrix describing the coupling structure, and $h: \mathbb{R}^d \to \mathbb{R}^d$ is the coupling function. We consider the individual node dynamics given by
HR \cite{Hidmarsh_Rose1984,huerta1998clusters} oscillators. For this model, $\mathbf{x}_i = 
\begin{bmatrix}
x_i,
y_i,
z_i
\end{bmatrix}^T$, and the individual node dynamics is
\begin{equation}\label{eqHR1}
f({\mathbf x}_i) = 
\begin{bmatrix}
y_i-ax_i^3+bx_i^2-z_i+I \\
c-dx_i^2-y_i \\
r(s(x_i-x_R)-z_i)
\end{bmatrix}.
\end{equation}
\ed{Above $x$ represents the membrane potential, $y$ is the rate of transfer of sodium and potassium ions through the fast channels, and $z$ is the adaptation current which reduces the spiking rate after a spike has occurred, see SI (Sec.~5.1) for more details about the parameters.} We consider diffusive coupling through all variables
\begin{equation}\label{diffusivecouple}
h_1(\bf{x_i},\bf{x_j})= 
\begin{bmatrix}
x_j-x_i \\
y_j -y_i \\
z_j-z_i
\end{bmatrix}.
\end{equation}
The diffusive coupling mimics electrical interactions between the neurons: a higher difference of '+' and '-' ions between pre-synaptic and post-synaptic neurons causes a proportionally higher flow of these ions through channels. 
We also consider a more realistic model of the neuronal dynamics, which includes coupling through two terms,
\begin{equation}
\label{eq:Model2b}
h_2(\mathbf{x_i},\mathbf{x_j})= 
\begin{bmatrix}
0 \\
y_j-y_i \\
0
\end{bmatrix}- \alpha(x_i-V_{syn})
\begin{bmatrix}
[1+e^{-\lambda(x_j-\theta_{syn})}]^{-1} \\
0 \\
0
\end{bmatrix}.
\end{equation}
The first coupling term in  Eq.~(\ref{eq:Model2b}) describes simple diffusive coupling through the $y$-variables only, while the second represents a ``chemical coupling" function.  This coupling scenario was presented in \cite{hizanidis2016chimera} as a more realistic consideration of two types of neuronal connections, one set which interacts through electrical signals and the other does so chemically. \paul{An interesting feature of this model was the coexistence of multiple different chimera states, even though the network did not} contain any non-trivial automorphism (symmetry) groups. Recently it has been shown that such symmetries are a sufficient \cite{Pecora2014,nishikawa2016network}, but not necessary \cite{schaub2016graph,gambuzza2019criterion} condition for a graph to support a stable chimera state. This is an important distinction since, in fact, the DTI network that we examine here contains no such non-trivial automorphism group. Indeed, as the number of nodes in a network increases, the lower the likelihood that the network will contain such symmetries \cite{kotters2009almost}. 

\noindent
\ed{\textbf{Simplification is the first step.} Our model of the brain dynamics incorporates simplifications, where we employ a single-neuron model to represent the dynamics of a node. While more complicated approaches such as the Wilson-Cowan nonlinear oscillator \cite{wilson1972excitatory,Bansal_2019} or the neural mass model \cite{Deco2015} \paul{could  better represent large pools of neurons}, the intricacies involved, such as higher-dimensional descriptions and noise, might obscure the essence of our observations. Addressing these challenges in more elaborate models is a task for future research.}

\subsection{Vector Pattern State}

\begin{figure}[t]
\includegraphics[width=1.00\linewidth]{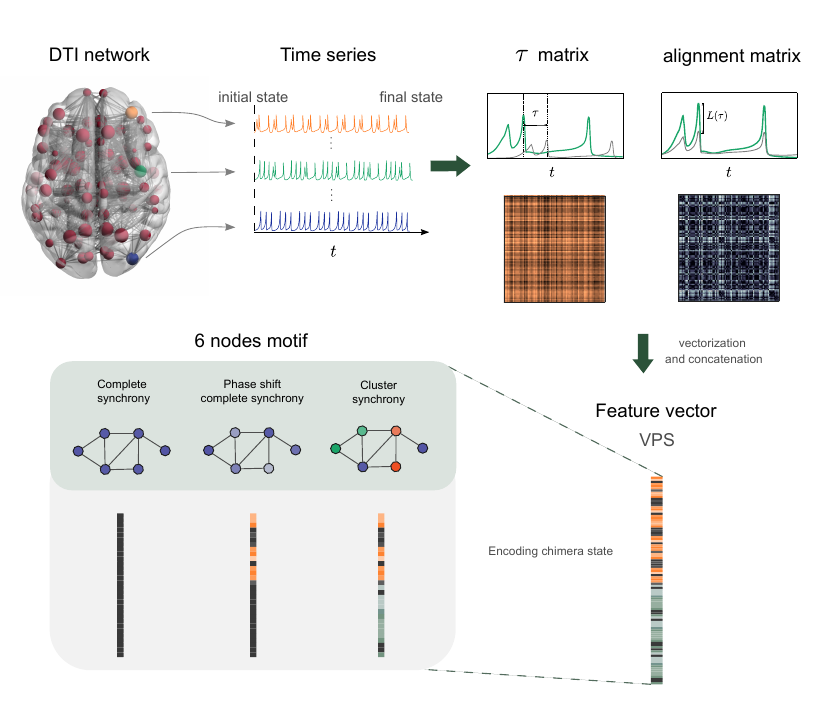}
\caption{\textbf{Schematic diagram of the Vector Pattern State construction.} (Top) \paul{The actual DTI network used in this work mapped to brain space, generated by BrainNet Viewer 1.7 (\url{www.nitrc.org/projects/bnv/}) \cite{Xia_2013}, is shown on the left. Nodes are structural brain regions and the edges are anatomical connections via white matter fiber tracts. The size of each node is scaled by the degree centrality.} From some initial state the dynamics of the three individual brain regions are shown as hypothetical time series, reaching a final state. The time shift $\tau$ and alignment between states of all pairs of nodes is recorded at the final state, yielding the $\tau$ and alignment matrix $L(\tau)$. (Bottom) To create a feature vector associated with this final state, we stack and concatenate these matrices into a single vector, defining the Vector Pattern State (VPS). The VPS encodes patterns of synchrony, with or without phase shift. All states correspond to different VPSs, and are here distinguished in the 6 node network, shown as different colored patterns.}
\label{StoryBoard}
\end{figure}

At some chosen initial time ($t=0$) the \ed{network} is in a particular initial state, see Fig. \ref{StoryBoard}. Each \ed{node} undergoes some dynamics, shown as a time series, and after a transient time, reaches a final state. Out of all time series generated by the \ed{network}, three are depicted in Fig. \ref{StoryBoard}. Each of the \ed{nodes} can be classified based on their level of activity by assigning each \ed{node} a color based on intensity, and \ed{nodes} with approximately the same level of activity are given the same color. 

A chimera state generally describes a scenario amongst $N$ coupled dynamical nodes \cite{abrams2004chimera,kemeth2016classification} whereby their time variables ${\mathbf z}(t)=({\mathbf x}_1(t),{\mathbf x}_2(t),...,{\mathbf x}_N(t))$ (in the notation here, ${\mathbf x}_i(t)\in {\mathbb R}^3$ denotes one of the coupled HR oscillators; in Eqs.~(\ref{generalcouple})-(\ref{eqHR1}), ${\mathbf z}(t)\in {\mathbb R}^{3N}$ encompasses the set of all the coupled variables) eventually converge to a state where some of the variables at nodes synchronize, $t>0$, possibly including a phase shift, while others of the variables are incoherent to those, but possibly synchronous amongst themselves. The latter scenario, with the remaining variables being synchronous amongst themselves, is also called cluster synchrony \cite{belykh2008cluster,Pecora2014}. 

\ed{Traditionally, activity patterns have been identified in terms of the level of synchrony of the overall system \cite{Deco2015,Bassett2017}. However, the system may exhibit synchronous, asynchronous, and partial synchrony, which encompasses chimera states. However, partial synchrony limits a richer characterization of the possible activity patterns. Indeed,} for a large system such as the DTI network of $N = 83$, the chimera states can be plausibly quite complex, with exponentially many plausible groupings, and many in fact are feasible. Thus, the characterization of different chimera states requires deciding which variables synchronize in the complex networked system of HR oscillators. 
 
To characterize a chimera state of the 83 brain regions, we quantify the level of synchrony between pairs of nodes in the network. More precisely, after a large time {\color{black}$T_0 > 0$} to allow transients to settle, the time series $x_i(t)$ are compared to $x_j(t-\tau)$ for each $i,j$ pair, as depicted in Fig. \ref{StoryBoard}. Allowing for phase shift synchrony by a possible shift, we must decide if 
\begin{equation}\label{Lfun}
L(i,j,\tau)= \lim_{T \to \infty} \frac{1}{T} \int_{T_0}^{T_0+T} \|{ \mathbf x}_i(s)- {\mathbf x}_j(s-\tau) \|_2^2 ds,
\end{equation}
is small for any phase shift $\tau>0$, which may be decided by minimizing $L(i,j,\tau)$. {\color{black} Here the limit to infinity means large enough integration time, see SI for practical implementation for finite time series.} Since the maximum of the cross-correlation has the property that,
\begin{equation}
    \underset{\tau}{\mbox{argmax} }({ \mathbf x}_i\star {\mathbf  x}_j)(\tau)=\underset{\tau}{\mbox{argmin }} L(i,j,\tau), \mbox{ each } i,j=1,2,\dots,N,
\end{equation} it is convenient to estimate when variables ${\mathbf x}_i(t)$ and ${\mathbf x}_j(t)$ settle into a synchronous state by maximization of the discrete cross-correlation,  
\begin{equation}\label{crosscor}
R_{x_i,x_j}(\tau) = \sum \limits_{t} x_i(t)x_j(t-\tau),
\end{equation}
in terms of the scalar $x_i$, the first index of each ${\mathbf x}_i$. 

After all pairs are taken into account, we construct the corresponding $\tau$ matrix and the {\color{black}alignment} matrix via $L(\tau)$. From these matrices, we create the feature vector, the vectorization and concatenation of the two matrices into a single vector, which we call the vector patterns state (VPS)
\begin{equation}
\label{VPS}
e_l=(\tau^*_{1,2},\tau^*_{1,3},\dots,\tau^*_{N-1,N},\beta L(1,2,\tau^*_{1,2}),\beta L(1,3,\tau^*_{1,3}), \dots ,\beta L(N-1,N,\tau^*_{N-1,N})),
\end{equation}
where the parameter $\beta\geq 0$ scales the importance of contrasting the optimal phase shift $\tau^*_{i,j}$ for comparison of the coupled components, and that best matched difference between components $L(i,j,\tau^*_{i,j})$. Whether complete synchrony, cluster synchrony, or chimera, with or without phase shift, all patterns are encoded via the VPS, as illustrated in Fig. \ref{StoryBoard}.

\subsection{\ed{Fractal basin structure \paul{supports} the nimble brain}}

\ed{Basin of attraction is defined as the set of all the initial conditions in the phase space whose trajectories eventually fall into a particular attracting state. In our case,} different initial conditions may lead to the same final state (and \ed{are assigned to the same color when visualized}) according to the VPS. It is the pairing of the initial state with the final state which we are interested in. This represents the structure of the basin of attraction to various final states. 

\ed{Recently, there has been significant research into unraveling the basin structure of attractors in high-dimensional systems \cite{Strogatz_2006,Menck2013,martens2016basins,Zhang_2021}. Typical questions about basin structure have centered around the size and shape of these basins, \paul{both quite} challenging in our specific case. We are dealing with a system comprising 83 nodes, each associated with a three-dimensional dynamical model, with a phase space that is $3\times 83=249$ dimensional. In contrast to many current studies that rely on characterizing states based on identical synchronization, our focus is on achieving approximate synchrony. We find this approach more versatile and applicable to a broader range of \paul{neuroscience questions where identical synchrony is unlikely.}} Hence, mapping the basin of attraction structure of the various chimera states \ed{based on approximate synchrony} becomes a problem of associating many long-time patterns from distinct initial conditions, and so this requires a way to match similar signals corresponding to occurrences of disparate chimera states.  The full basin structure is too complex to visualize, hindering any chance to uncover \ed{its structure, and consequently,} the mechanism of the nimble brain. To this end, we use the introduced VPS to solve this mapping problem.

We wish to partition a randomly selected ``slice" of the phase space into those regions with similar asymptotic behavior, by observing a sample of $M$ initial conditions which we index by $l$, ${\mathcal Z}=\{{\mathbf z}_l(0)\}_{l=1}^M$.  To this end, we wish to decide the synchrony pattern of any one ${\mathbf z}_l(0)$, by comparing the long time state of component time series according to Eq.~(\ref{Lfun}) at optimally matched phase shift, according to Eq.~(\ref{crosscor}). With the VPS, we can now assert that two initial conditions ${\mathbf z}_{k_1}(0)$ and ${\mathbf z}_{k_2}(0)$ yield asymptotically similar complex synchrony patterns only if their VPS are relatively close, i.e. $\|e_{k_1}-e_{k_2}\|_2$ is small. 

\begin{figure}[t]
\includegraphics[width=1.00\linewidth]{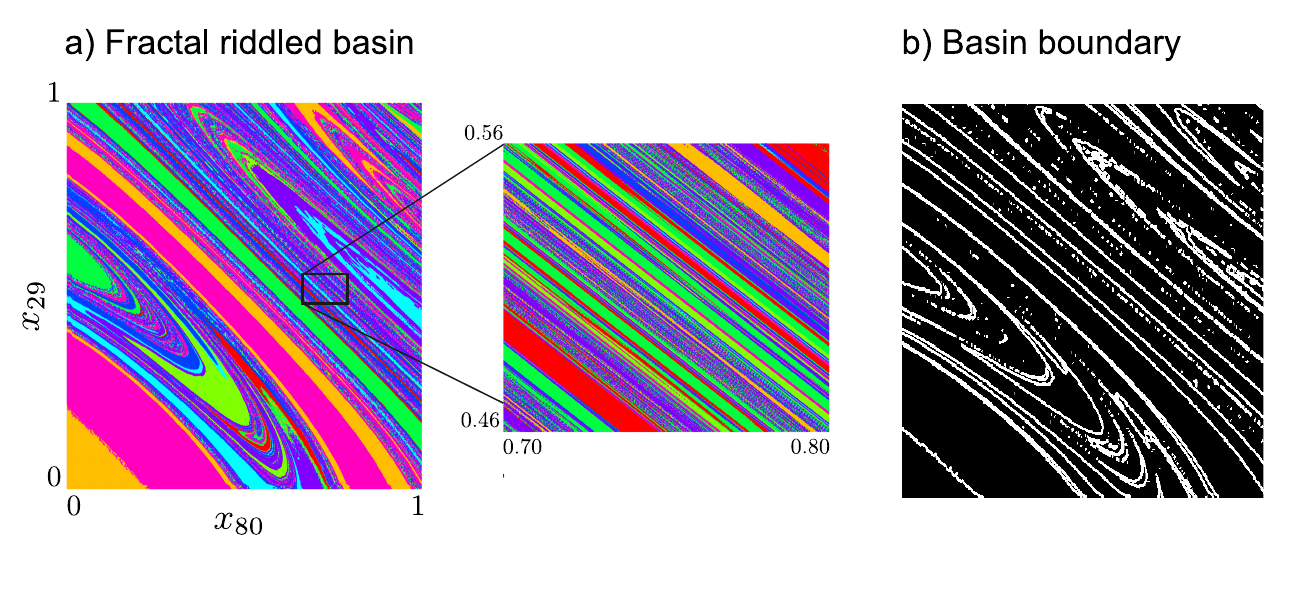}
\caption{{\color{black}\textbf{Fractal riddled basin of the full featured HR oscillator model on the DTI network.} a) An arbitrary plane ``slicing" through the full high dimensional space was selected on which initial conditions are sampled uniformly. Here the $x$ component of the $29^{th}$ oscillator and the $x$-component of the $80^{th}$ oscillator, at $t = 0$ define the plane. In this basin, the initial conditions associated with different chimera are each a different color. Note that in a region that appears to alternate between just a few states, actually exhibits a rich structure with many different interleaved states when zoomed in at higher resolution. b) The basin boundary set shown in a). The box counting fractal dimension of the basin boundary in this plane, which is computed $d_B\sim 1.8$, being non-integer indicates a fractal set. We consider full featured HR oscillator model Eqs.~(\ref {eqHR1}),(\ref{eq:Model2b}) with $a=1,b=3,c=1,d=5,s=4,r=0.005,x_R=-1.6,I=3.25,\sigma=0.5,\alpha=0.03, V_{syn} = 2, \theta_{syn} = -0.25$ and $\lambda = 10$. The partition into basin structure associated with distinct dynamical chimera states follows k-means clustering on the VPS structure, Eq.~(\ref{VPS}), using the cost Eq.~(\ref{Lfun}), inferred with cross-correlation, Eq.~(\ref{crosscor}), using $k=8$, the result of a classic elbow method.}}
\label{fractalfig}
\end{figure}

Now the problem of partitioning the phase space into like asymptotic chimera states reduces to a clustering problem of all VPSs relative to the different initial conditions.  To this end we apply the k-means method to the set of VPS,  $\{e_l\}_{l=1}^M$, to cluster the space into k-regions (colors) and we map the phase space by associating these colors to each corresponding initial condition ${\mathbf z}_l(0)$.  Thus the clustering is a partition function, ${\cal P}:{\cal Z}\rightarrow \{1,2,..,k\}$, as shown in Fig.~\ref{fractalfig}.  We describe these as basin plots since in any like colored region, the orbits of the initial conditions map asymptotically to similar patterns.  
Relevant details concerning the experimental methods are included in the figure caption. \ed{ As noted above, a key component of our method in determining how to group the final states into their various attractors is clustering. While numerous clustering methods exist, we chose, for reasons of computational complexity, k-means. Thus a} general description of the k-means algorithm as a clustering method, and the manner in which we choose how fine to partition the space with the selection of a specific $k$ are both presented in the SI.

%\reformulate{Jeremie: Maybe further expand or comment that details are in the Methods section. So, make it easier to see the influence of k-means for identifying the basin. }

\subsection{Coupled HR oscillators in a DTI network}
Even with these simplified dynamical models of the brain, there is still rich complexity that demonstrates interesting phenomena in the basin structure. In \ed{Fig.} \ref{fractalfig} we show that using the coupled HR oscillator model, the basin boundary between the states has a non-integer Hausdorff dimension, and thus fractal basin boundaries. \ed{In the parameter regime $a = 1, b = 3,c = 1,d = 5,s = 4,r = 0.005,x_R = -1.6,I = 3.25,\sigma = 0.5,\alpha = 0.03, V_{syn} = 2, \theta_{syn} = -0.25$ and $\lambda = 10$, which is known to contain chimeras \cite{hizanidis2016chimera}, we use the electrical and chemical coupling functions, Eq. \eqref{eq:Model2b}, where the corresponding adjacency matrices are assumed to be the same, unlike in \cite{hizanidis2016chimera}.} Here for the first time, we map the manner in which these states are intricately co-mingled.  On an arbitrary plane, in this case, which we selected randomly as a slice of the full phase space restriction for the sake of visualization, a uniform grid of $750 \times 750$ initial conditions is chosen. The various colors label initial conditions associated with differing chimera state states. Furthermore, ``zoom" restrictions of the domain are also shown to illustrate the fractal-like structure of the basins of attraction at a finer scale. \ed{We validate this assertion by computation, that the basin boundaries projected into the planes shown to have a box counting dimension that is not an integer. The box counting dimension of the boundary sets was found to be fractal in Fig. \ref{fractalfig} (b), where the dimension was estimated to be $d_{\mbox{box}}\sim 1.8,$ by the method described in Eq.~\ref{eq_box_dimension}.}  

The basin structure in Fig. \ref{fractalfig} appears to exhibit complexity beyond simple fractal basin boundaries. A riddled basin structure appears, which is the scenario that regions exist where points in the domain of one attractor have the property such that small neighborhoods of nearby points have a nonzero probability of being in the basin of another attractor \cite{alexander1992riddled,ott1994transition,cazelles2001dynamics}. In practical terms, this means that there are large regions in phase space where it is likely that even small perturbations can send the outcome to regions corresponding to a different state.  This has significant implications for the possibility of nimble switching between states, since switching between multiple states that may be co-mingled in the phase space may require only vanishingly small control inputs.

%\reformulate{Jeremie: There are two aspects to discuss here based on Reviewer 2 comment:
%\begin{itemize}
 %   \item Eq(3) defines individual node dynamics as electrical synapses, but Eq(4) introduces a more realistic model of synapses that needs clarification regarding the use of either Eq(3) (Fig3) or Eq(4)(Fig2) in computations.
 %   \item  Fig.2 and 3. Could authors explain the numerical value associated with the fundamental set of parameters?
%\end{itemize}
%}

\ed{\subsubsection{Fractal basins are ubiquitous}}

\begin{figure}[t]
\includegraphics[height=\linewidth, width=\linewidth]{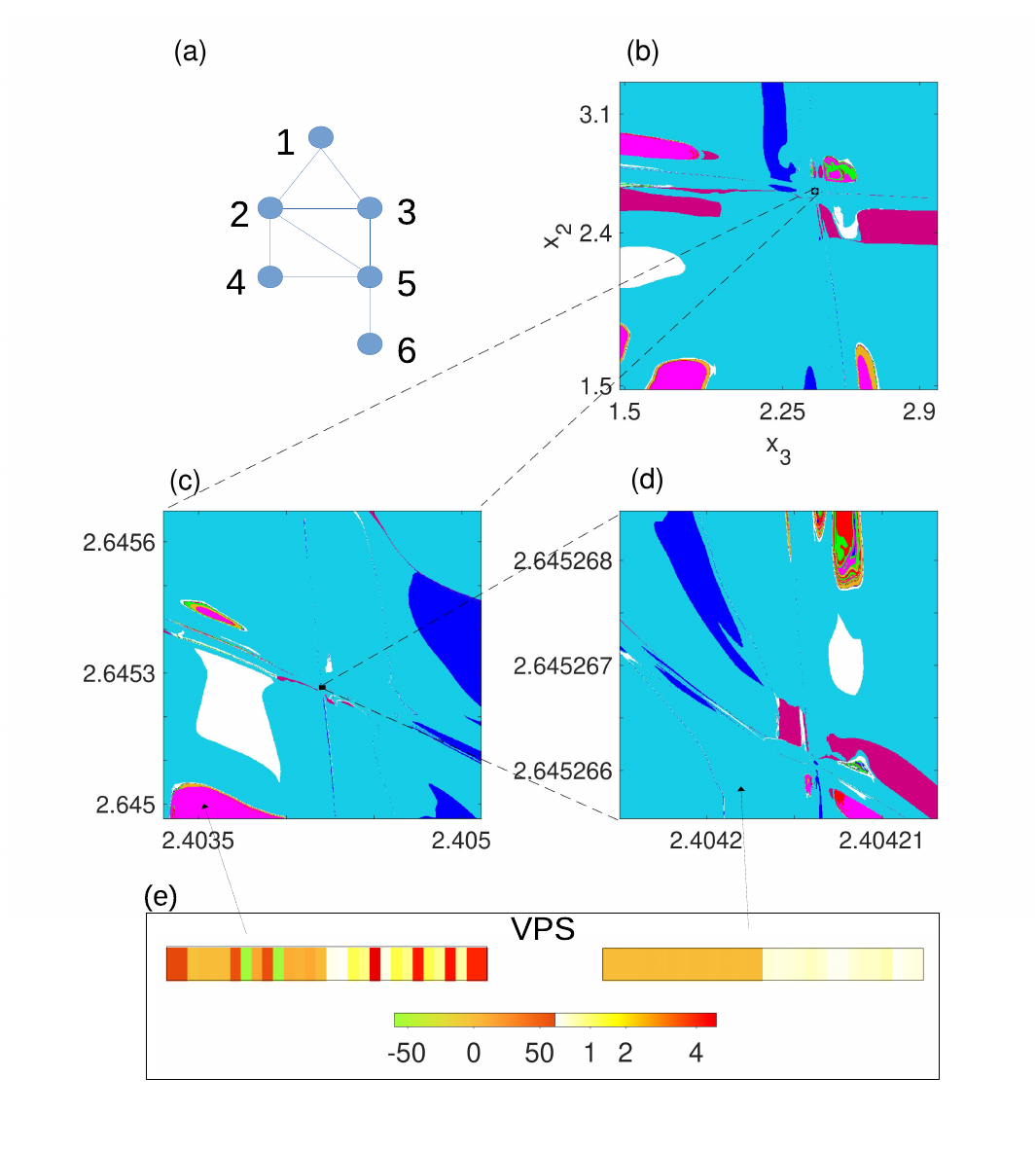}
\caption{
A simplified HR model with diffusive coupling Eqs.~(\ref{eqHR1})-(\ref{diffusivecouple}) on a small graph illustrates the ubiquity of fractal basin structure of chimera states.
(a) A network of $6$ nodes that does not contain non-trivial symmetries. 
Nonetheless, there are many stable chimera states (at least on the time scale examined), and the basin structure shown in $8$ colors indicates distinct patterns that can be derived by VPS structure, Eq.~(\ref{VPS}), by the same method as in Fig.~\ref{fractalfig}. (b) Fractal basins for HR oscillators on this network when $x_R=-0.5(1+\sqrt{5}), I=3.27, r=0.017, \sigma =0.0004$, and {\color{black}$\beta=1$}. All other $x_i, y_i$, and $z_i$ values at $t=0$ are initialized to be $-0.5$. (c) and (d) are zoomed regions indicated by the black rectangles in (b) and (c). \anil{(e) Centroid locations of two of the clusters in $\tau-L$ space, which resembles the approximate form of most of (or all) VPSs inside (see SI for a detailed view of all $e_l$ vectors inside each cluster).}}
\label{fig_2}
\end{figure}

\ed{
\noindent
\textbf{HR oscillators coupled in small networks.}} To illustrate the generality of our results, \ed{we present fractal basins in different networks. } Fig. \ref{fig_2} displays complex patterns that can be found in the basin of a smaller network of 6 oscillators, as shown in Fig. \ref{fig_2} (a). \ed{We use the electrical coupling scheme with $h_1$ given in Eq.~\eqref{diffusivecouple}, and the parameter values based on earlier research works, see \cite{Djeundam2013, Hindmarsh1984}.} We chose to examine a small synthetic network, which does not have any non-trivial automorphism group, to demonstrate the ability of a coupled HR model to form a basin that has fractal boundaries. \ed{In fact, in Fig.~\ref{fig_2} (b), the corresponding estimate is $d_{\mbox{box}}\sim 1.27$, where it} shows the basin structure grouped into 8 different states using k-means. Figs. \ref{fig_2} (c) and (d) are shown in zoomed (restricted) in regions of Fig. \ref{fig_2} (b) and Fig. \ref{fig_2} (c). The structure of the basin is quite complex at all scales examined. 
% In these cases there appear to be fractal basin boundaries. 

\begin{figure}[t]
\includegraphics[height=\linewidth, width=\linewidth]{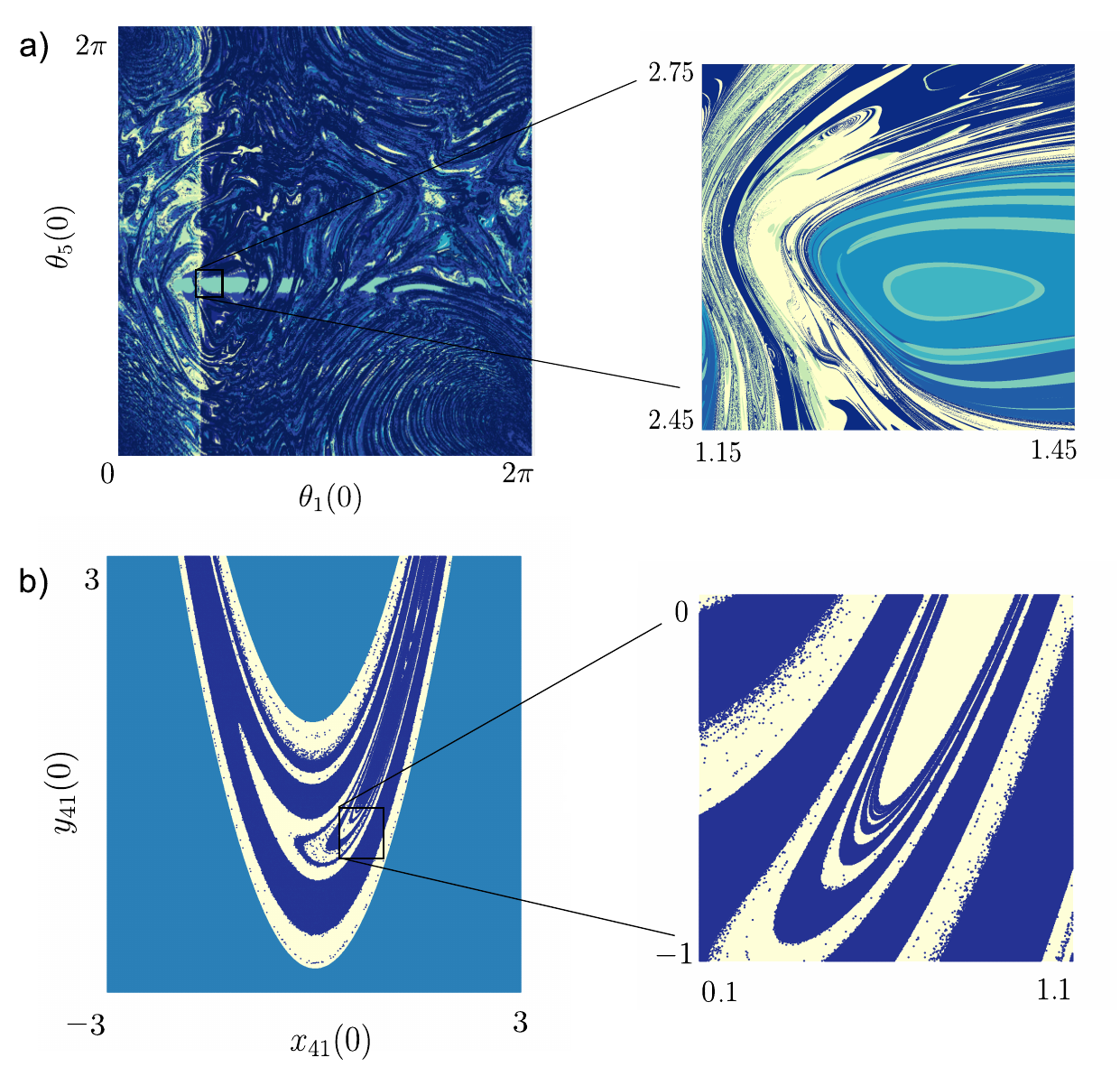}
\caption{{\color{black} Riddled basins for different networked systems. (a) The left panel shows a two-dimensional section of the state space for a system of coupled phase oscillators on a network showing basins of 12 (clustered) distinct states. Right panel Zoomed in from inset of a)  showing basins of 7 (clustered) distinct states. To construct the VPS, we use $\beta = 1$ in Equation \eqref{VPS} and a grid with $1248 \times 1248$ and $624 \times 624$ for left and right panels, respectively, uniformly sampled initial conditions. (b) H\'{e}non map dynamics on a DTI network with no non-trivial symmetry. See further details in the SI.} }
\label{fig:Figure4}
\end{figure}

We further explore two more examples of local dynamics and network structure to support the generality of our claims on the nimble brain. In Fig.~\ref{fig:Figure4} we illustrate these examples, and thus the ubiquity of complex basin structure between various chimera states.  

\ed{
\noindent
\textbf{Identical Kuramoto oscillators.} We consider the following equations of motion for the identical oscillators
\begin{align}
\dot{\theta}_{i} = \sigma \sum_{j = 1}^{N} [A]_{i,j} \sin(\theta_j - \theta_i - \alpha), \quad i = 1, \dots, N, 
\end{align}
where $\sigma$ is the overall coupling strength and $\alpha = \pi/2 - \ed{\gamma}$ with $\ed{\gamma} = 0.025$. The adjacency matrix $A$ represents a network that does not have full permutation symmetry. To generate this network we initiate two populations of 5 nodes that are globally coupled akin to \cite{montbrio2004synchronization}, and remove uniformly at random one edge from the graph, see details in the SI. Fig.~\ref{fig:Figure4}(a) shows the complex basin structure that is captured using our VPS.

\noindent
\textbf{H\'{e}non map.} Additionally, we study the network of coupled H\'{e}non maps,
\begin{equation}
\begin{bmatrix}
x_{i}({t+1}) \\
y_{i}({t+1})
\end{bmatrix} = 
\begin{bmatrix}
f_x(x_{i}(t),y_{i}(t)) + \sigma \sum \limits_{j = 1}^N[A]_{i,j}\Big(f_x(x_{j}(t),y_{j}(t)) - f_x(x_{i}(t),y_{i}(t))\Big) \\
f_y(x_{i}(t),y_{i}(t))
\label{eq:HenonMap}
\end{bmatrix}
\end{equation}
for $i \in \{1,2,..., N\}$, with $f_x(x,y) = 1-px^2+y, f_y(x,y) = bx$ and $t \in \mathbb{N} $, as discussed in \cite{santos2018riddling}. The parameters chosen are $p = 1.44$, $b = 0.164$, $\sigma = 0.8$. The network used is the DTI brain network from Fig. \ref{StoryBoard}. Fig.~\ref{fig:Figure4}(b) again highlights the generality of the complex structures and also the utility of the VPS technology.   Further details of both of these examples are presented in SI (Secs. 5 and 6).
}

\section{Discussion}

The brain has proven to be extremely nimble in its ability to switch between states in response to stimuli, thoughts, and/or decisions. As observed by various imaging techniques, this is associated with rapid switching between patterns of synchronous, chimera, and incoherent states.

\noindent
\ed{\textbf{Basin structure of network dynamics}.} Several prior works have studied the basin structure of chimera states in networked systems. There have been observations of chimera states with an intermingled basin structure in a special case of a strongly self-coupled cluster network specifically designed to emphasize chimera; see an explanation of critical switching behavior \cite{zhang2020critical}. Authors in \cite{brezetsky2021chimera} found highly riddled basins in small and highly symmetric all-to-all networks of coupled phase oscillators. \ed{Fractal basins of chimeras states were found in small networks of coupled complex maps \cite{Andrzejak_2021}.} In \cite{martens2016basins} the authors use a low-dimensional description valid for the infinite size system \cite{ott2008low} to characterize the basin structure of different patterns in a model of two populations of all-to-all coupled Kuramoto oscillators \cite{montbrio2004synchronization}. Likewise and related, in \cite{li2022basins} analyze the same highly symmetric two population network model for chimera, but then illustrate chimera states for a DTI network with coupled Wilson-Cowan oscillators. They define chimera states in terms of a highly approximate synchrony, which is not a general approach such as our VPS that would allow for analysis of basin structure. Similarly, in \cite{santos2017chimera} chimera premised on approximate synchrony was described for a cat brain connectome data set \cite{scannell1993connectional} describing coupled HR oscillators as coupled through one variable only, but again, no basin structure was found. In \cite{santos2018riddling}, authors use the chaotic H\'{e}non map coupled by again a highly symmetric network, the circulant (ring) stricture, and thus to find fractal basins for chimera premised on identical synchrony. 

%\ed{  \textcolor{violet}{(Paul says- I am a bit confused by this next sentence. Reference 22 used DTI and WC oscillators. that seems to do exactly this, but they cant get at the basin structure.)} However, none of these works is capable of handling a fully asymmetric system as is naturally presented by a biologically true brain, such as a human DTI derived network, with realistic dynamics such as we do with the HR model.}

% Chimera patterns also on a DTI derived network, but specifically on FitzHugh-Nagumo neurons instead, was studied in \cite{chouzouris2018chimera} for the sake of understanding epileptology and epilepsy, but the present work is the first to take a broader view of healthy function in terms of completely mapping basins of stability of various pattern states.

\noindent
\ed{\textbf{Dynamical systems theory is useful to explain the brain.}} \ed{Dynamical systems theory has been adopted as an approach to gain insights over the brain dynamics across various scales \cite{BRAUN2010740,Breakspear2017,esteban2018informational,Graben_2019,morrison2022chaotic,Yohan_2022,SCHIRNER_2022,Tsuda2022}. Instead of an empirical or quantitative investigation, e.g. trying to observe attractor-like states \cite{BRAUN2010740,Balaguer_2011}, most investigations have focused on proposing theoretical dynamical mechanisms \cite{Graben_2019,Yohan_2022}. For example, dynamical systems theory has contributed to the development of } theories of consciousness, by so-called integrated information theory (IIT) \cite{esteban2018informational}, or the description of complex switching phenomenon in biological systems by the concept of chaotic heteroclinicity \cite{ morrison2022chaotic}. 

\ed{Within a dynamical systems perspective, numerous possible mechanisms exist, necessitating research to pinpoint the one that aligns most closely with empirical data. In this context,} we provide numerical evidence of fractal basin boundaries that have non-integer box counting dimension, and riddled basin boundaries.\ed{This evidence corroborates a theoretical explanation for resting-state brain dynamics, as investigated in \cite{Graben_2019}, which shows the promise of this dynamical mechanism. } We observe these properties in numerical simulations of \paul{multiple different systems of coupled \ed{dynamical}} oscillators, using an experimentally determined \paul{human structural brain network as well as} small test networks. With this evidence, we have identified a  \paul{potential mechanism that would allow a nimble brain to switch between various distinct states with only small changes in the system parameters.}

\ed{From a dynamical systems perspective, we argue that} coexisting attractors corresponding to the various chimera states may seemingly suggest that large perturbations would be required to transition from deep in the well of one stable state to another. \paul{A brain with such dynamics} would be at odds with the idea of a system that can nimbly switch between states. From a neuroscience perspective, it may seem that to transition from one brain state to a distinctly different brain state, one would have to traverse many unique states on a trajectory to the final desired state. We offer an explanation for how to resolve this seeming contradiction in the form of fractal basin boundaries. The fractal basin boundary allows for different stable states to be mixed together closely, creating the opportunity for small perturbations to lead to entirely different stable states, as patterns of chimera.  

Thus, the main results of this work are summarized as follows:
\begin{enumerate}
\item Our main proposal is that brain activity switching, that is, the nimble brain, is explained by fractal intermingled (riddled) basins. Complex basins of attraction for each chimera state are intrinsically highly intermingled. Thus, significantly different states are nonetheless near each other, in the dynamical variables of the phase space, and so available for nimble control manipulations by internal cognitive processes or external environmental events.  
\item Even though the networks in the system have no symmetries, a generalized interpretation of synchrony allows fractal (intermingled) riddled basins, including relatively small model networks.
\item A crucial technology that underpins these above two assertions is based on clustering the VPSs corresponding to chimera states. Here, the k-means of a metric between VPS is a convenient clustering approach. Implementation of the computational task in mapping fractal basins is a key technical innovation that we have developed as background for this new description of the neuronal dynamics of the brain. \ed{Our approach can be extended to more complex models of brain dynamics.}
\end{enumerate}

\ed{Our approach allows the first step to find basin structure of complex high-dimensional systems. \paul{Our initial description of such fractal basins necessitated a somewhat simplistic, though biologically inspired, brain model. Now that we have presented this potential mechanism for nimble brain state shifts, experimental neuroscientific studies are needed to empirically validate, or reject, the hypothesis that we have presented. We also envision studies that further investigate the structure of these basins.} Promising directions include octopus-like basins for basin structures for chimera states \cite{Zhang_2021}, narrowing down other potential mechanisms for the nimble brain.}

\section{Methods}
\label{sec11}

%\subsection{\ed{Clustering methods to construct VPS}}

%\reformulate{Jeremie: I would suggest we include some details of the clustering methods to satisfy Reviewer 2. For instance, Reviewer 2 could not find the details of the k-means. I would suggest to bring large part of the SI to here, in particular the discussion that other methods could be employed. Our proposal was to use the simplest as possible to illustrate our findings. }

\subsection{Fractal basins: box counting dimension}
The assertion of fractal basin boundaries is a matter of considering the approximate boundary set $S_{BL}$, such as the one shown in Fig.~\ref{fractalfig}(b), from the basin set in Fig.~\ref{fractalfig}(a), shown in cross-section with respect to the variables.

The box counting dimension can be estimated by counting a covering of squares of side length $\epsilon$, and then consideration of this count $N(\epsilon)$ upon refinement by decreasing $\epsilon$.
The box dimension is defined \cite{dubuc1989evaluating}:
\begin{equation}
    \mbox{d}_{\mbox{box}}(S_{BL}) = \lim \limits_{\epsilon \rightarrow 0}\frac{\mbox{ln}(N(\epsilon))}{\mbox{ln}(1/\epsilon)},
    \label{eq_box_dimension}
\end{equation}
that is equivalent to the Minkowski-Bouligand dimension.
While $S_{BL}$ is simply a slice of the full high-dimensional boundary set, the non-integer result, $\mbox{d}_{\mbox{box}}(S_{BL})=1.8$, together with the statistically self-similar structure shown, supports the assertion of a fractal set. {\color{black}  Likewise, in Fig.~\ref{fig_2}(b), the corresponding estimate is $d_{\mbox{box}}\sim 1.27$. }

\backmatter

\bmhead{Data availability} 
The network structure used here was derived from diffusion tensor imaging, and parcellated by the Lausanne anatomical atlas into 83 anatomical regions. 
This structure is publicly available  \cite{bonilha2015reproducibility}, see the link \url{https://rb.gy/q3o71}, from which we selected ``Subject 1" as used in \cite{fish2021entropic}. \ed{The visualization of the DTI network is generated by BrainNet Viewer 1.7 (\url{www.nitrc.org/projects/bnv/}) \cite{Xia_2013}.}

\bmhead{Supplementary information}

A supplementary information file of further details, theory, and explanations is included.

\bmhead{Author contributions} 

E.B., J.F., A.K. and P.L. designed research, and all authors performed research and analyzed results. J.F., A.K. and E.R.S. implemented the numerical simulations. All authors contributed to the writing of the manuscript. All authors reviewed and approved the final manuscript.

\bmhead{Competing interests}
The authors declare no competing financial interests.

\bmhead{Acknowledgments}

A.K, J.F., E.R.S., P.L. and E.B. acknowledge support from the NIH-CRCNS. J.F. and E.B. are also supported by DARPA RSDN. Additionally E.B. is supported by the ONR, ARO and AFSOR. E.R.S. was also supported by Serrapilheira Institute (Grant No.
Serra-1709-16124). 
% The brain graphic in Fig.\ref{StoryBoard} was extracted from FreePick.


\bibliography{sn-bibliography}% common bib file

\newpage


\section*{Supplementary material}

# Fractal Basins as a Mechanism for the Nimble Brain

Erik Bollt, Jeremie Fish, Anil Kumar,
Edmilson Roque dos Santos and Paul Laurienti

October 8, 2023

# Contents

Here we expand concepts with further details as they appear in the body of this work. Specifically, network concepts for complex systems, dynamical systems, fractals concepts, and machine learning concepts for neuro-inspired models of the brain.

# 1 Networks and Complex Systems

## 1.1 Networks

Networks play an important role in many complex systems, including the brain. A networked representation of the human brain, where different regions are treated as vertices and interactions between them as edges, allows us to explore both the structural and dynamical properties.

In this work, we use the phrase network synonymously with the mathematical conception of a graph $G = (V(G), E(G))$ which is a set of vertices (nodes) $V = \{v_1, v_2, ..., v_N\}$, with cardinality $|V(G)| = N$ and edges (links) $E \subseteq V \times V$ with $|E(G)| = M$. The adjacency matrix $A$ is a matrix representation of the graph with entries

$$[A]_{i,j} = \begin{cases} w_{i,j} \text{ if } (v_i, v_j) \in E \\ 0 \text{ otherwise,} \end{cases} \quad (1)$$

where $w_{i,j}$ is a weight. A simple graph is an unweighted graph with no self-loops or multiple edges. Thus, each entry $[A]_{i,j} \in \{0, 1\}$ and $[A]_{i,j} = 0$ if $i = j$.

## 1.2 Subgraphs

Informally, a subgraph is a graph that results from the erasure of some of the vertices and edges of the original graph. That is, a graph $H = (V(H), E(H))$ is called a subgraph of a graph $G$ if $V(H) \subseteq V(G)$ and $E(H) \subseteq E(G)$ and $|V(H)| = N' \leq N$, $|E(H)| = M' \leq M$. Now suppose the vertices in $G$ are relabeled so that the corresponding vertices in $H$ have the same label, that is $V(G) = \{v_1, v_2, ...v_N\}$, and $V(H) = \{v_1, v_2, ..., v_{N'}\}$. Then we call

$H$ an induced subgraph of $G$ if $E(H) = (V(H) \times V(H)) \cap E(G)$. In other words, if for every vertex in $H$, the edges formed between each pair of those vertices in the graph $G$ are also found in $H$, then $H$ is an induced subgraph. A graph is called a proper subgraph if it is a subgraph of $G$ and it is not identically the same as $G$, i.e., it has strictly fewer nodes than $G$. Suppose a graph is partitioned into two proper induced subgraphs $H, K \subset G$ such that $V(H) \cup V(K) = V(G)$ and $V(H) \cap V(K) = \emptyset$; that is, $H$ and $K$ are mutually exclusive induced subgraphs. Then $G$ is called connected if for every such $H$ and $K$ pair in $G$, $E(H) \cup E(K) \neq E(G)$. That is, if for every pair of proper induced subgraphs satisfying the properties above, there is at least one edge between the pair, then $G$ is connected. Said another way, there is a path from every vertex to every other vertex. If $G$ is not connected, then it is called disconnected.

## 1.3 Graph Symmetries

Symmetry is an important property of graphs that feature prominently in the topic of dynamics on networks, specifically synchrony, though it is not a necessary condition [42]. By symmetry, we mean an automorphism of a graph. A graph automorphism is a permutation of the vertices that preserves the neighbors of any given vertex, meaning the graph remains "essentially" the same. Stated in terms of the adjacency matrix $A$ of a graph, there is a similarity transformation by a permutation matrix $P$ such that $A = P^T A P$, where $P$ can be constructed from the identity matrix by permuting rows and columns only, also here $P^T$ is the transpose of $P$. Such an automorphism is called a symmetry of a graph if it exists and the graph is called symmetric. Note that the term symmetric here does not mean a symmetric adjacency matrix, but the phrase is usually only used when there exists such a non-trivial permutation. Also, $P \neq I$, the identity matrix, which otherwise always works.

## 1.4 Diffusion Tensor Imaging (DTI) Graphs

Diffusion Tensor Imaging (DTI) is a magnetic resonance imaging technique used to find the structural organization of white matter tracts in the human brain. Since its formation invention in 1994, the technique has been widely used in brain research. The DTI data can be processed to estimate the brain fiber tracts using a technique called tractography. For the graph shown in Fig.1 of the main text, which is Subject 1 from [9], a process called deterministic tractography, which is available in a software package [45] called the diffusion toolkit, was performed. The brain was parcellated into 83 regions of

interest (ROIs) and fiber tracts from each were computed. Using the ROIs as individual vertices and fiber tracts between ROIs as edges, a weighted structural brain network was generated. The DTI network is then converted to a simple graph as seen in Fig.1 of the main text by changing all non-zero weights to 1 as was done in [13]. When we refer to a DTI network throughout the main text and in the SI, it is this graph that we are referring to.

# 2 Dynamics On Networks, Synchronization, and Chimera States

## 2.1 Synchronization

In a coupled dynamical system, synchronization means a '*rhythm*' or '*a kind of sympathy*' between two or more of its dynamical entities. The simultaneous firing of neurons in the human brain, synchronized flashing of light by fireflies, and synchronized clapping by an audience are some of the popular examples [34]. Synchronization arises due to interactions between dynamical elements: while the individual units try to move according to their own set of rules, attractive interactions or coupling with others try to bring harmony between them, which in turn leads to synchronization. Below, different forms of synchronization are discussed, each of which serves a role in our analysis of patterns in the brain, as deduced by our Vector Pattern States (VPS) to distinguish basins.

## 2.2 Identical Synchronization

The most basic and perhaps standard version of synchronization is called "identical synchronization" (IS) or "complete synchronization" (CS). In a complex system of coupled dynamically evolving units, this state occurs when all units asymptotically "beat together", identically.

Coupled elements are each modeled as differential equations with a vector field defined by a function $f$ described on a graph as vertices, with the edges describing communication between the vertices, through a coupling function $h$. A generic model used for study in such systems in continuous time is written as:

$$\dot{\mathbf{x}}_i = f_i(\mathbf{x}_i) + \sigma \sum_{j=1}^{N} [A]_{i,j} h(\mathbf{x}_i, \mathbf{x}_j). \qquad (2)$$

Here, $\dot{\mathbf{x}}_i$ represents the derivative of $\mathbf{x}_i \in \mathbb{R}^d$ with respect to time, the pair $\dot{\mathbf{x}}_i = f_i(\mathbf{x}_i)$ represents the isolated dynamics associated with an individual

vertex, $[A]_{i,j}$ are entries of the adjacency matrix and $h : \mathbb{R}^d \to \mathbb{R}^d$ is the coupling function between a pair of vertices. In the case where we have $f_i = f$ ($\forall i$) (meaning that the uncoupled oscillators obey identical dynamics) and $h(\mathbf{x}_i, \mathbf{x}_j) = \mathbf{0}$ when $\mathbf{x}_i = \mathbf{x}_j$ a synchronization manifold can be defined as

$$\mathcal{M} = \{(\mathbf{x}_1, \mathbf{x}_2, ..., \mathbf{x}_n) | \mathbf{x}_1 = \mathbf{x}_2 = \cdots = \mathbf{x}_n\}. \quad (3)$$

Any collective state of all units which occurs in $\mathcal{M}$ is called identically synchronous. Generally, states that asymptotically converge to the synchronization manifold are of interest, initial conditions which do so are defined as being in the basin of synchronization. Thus the set,

$$\mathcal{B} = \{(\mathbf{x}_1(0), \mathbf{x}_2(0), ..., \mathbf{x}_n(0)) | (\mathbf{x}_1(t), \mathbf{x}_2(t), ..., \mathbf{x}_n(t)) \in \mathcal{M} \text{ as } t \to \infty\}, \quad (4)$$

of all initial conditions that are asymptotically attracted to the synchronization manifold is known as the synchronization basin.

**Stability of identical synchronization**

In their seminal work relating to identical synchronization, Pecora and Carroll [33] analyzed the (local) stability of the synchronous state of a network coupled system. Let $\vec{x}(t) = (x_1(t), x_2(t), ..., x_n(t))$, then a synchronous state is called stable if starting with $\vec{x}(t) \in \mathcal{M}$, there exists a $\delta > 0$ and $\vec{x}(t) + \vec{\Delta} = (x_1(t) + \Delta_1, x_2(t) + \Delta_2, ..., x_n(t) + \Delta_n)$, with $\Delta_i \neq \Delta_j$ if $i \neq j$ and $|\Delta_i| < \delta$ such that $\vec{x}(t) + \vec{\Delta} \in \mathcal{B}$. That is the synchronous state is called stable if an infinitesimal perturbation to that state away from the synchronization manifold returns to $\mathcal{M}$ asymptotically. In [33] details of the conditions under which an identical system will have a stable synchronous state are stated.

## 2.3 Phase Synchronization

Phase synchronization between two signals $x_i(t)$ and $x_j(t)$ is a generalization of identical synchronization. Whereas for identical synchronization, $x_i(t) \to x_j(t)$, as $t \to \infty$, for phase synchronization, a slightly weaker but otherwise similar condition must hold after one of the signals is shifted by a phase. If there exists a phase shift, $\tau > 0$ such that $\|x_i(t) - x_j(t - \tau)\| \to 0$ as $t \to \infty$, then the two signals are defined as phase synchronous. This property is prominent in our analysis of brain patterns.

## 2.4 Generalized Synchronization

Generalized synchronization is another weakening of the concept of synchrony. Two signals $x_i(t)$ and $x_j(t)$ are defined as generalized synchronous if

$\|x_i(t) - f(x_j(t))\| \to 0$ for some nonidentity function $f(x) \neq x$. In the identity case that $f(x) = x$, then this would simply be a statement of identical synchronization. Usually, the statement is in terms of a smooth function $f$, and the set of points $\mathcal{D} = (x, f(x))$ on which orbits $(x_i(t), x_j(t))$ are attracted is called the synchronization manifold, and this phrase can be used even in the identical case. It is certainly possible for a system to exhibit phase and generalized synchrony.

## 2.5 Approximate Synchronization

If two signals $x_i(t)$ and $x_j(t)$ approach a synchronous state, such as identical synchrony, but they do not converge to that state, then the synchronous state is stable but not asymptotically stable. In other words, if for some time $T > 0$, and small $\epsilon > 0$, then $\|x_i(t) - x_j(t)\| < \epsilon$ for all $t > T$, but they do not converge to zero, then the signals are approximately (identically) synchronous.

## 2.6 Chimera States

A chimera state is defined by the presence of both synchronous and asynchronous behaviors in a complex system. Specifically, if a collection of coupled dynamical systems with corresponding signals, $\{x_1(t), x_2(t), ..., x_n(t)\}$ results in at least one subset of these that synchronizes (in terms of one of the types of synchronization noted above), and a distinct other subgroup that does not synchronize, then this is called a chimera state. Amongst those that synchronize, there may be just one large synchronous cluster, or perhaps several clusters that synchronize but not with each other, this is still called chimera, as long as there is also a subset that is asynchronous.

Chimera states have long fascinated the synchronization community, being first named chimera states (in an identical synchronization context) by Abrams and Strogatz in 2004 [2] though this phenomenon was reported earlier by Kuramoto and Battogtokh [20]. Though originally the term chimera state was reserved for identical oscillators, the literature is now filled with examples in which the term chimera state is used for non-identical oscillators [21, 38, 27, 18, 22] and so we have chosen to use the term in the broadest sense. The definition used here allows states that are sometimes called cluster synchronization or simply a clustered state [19, 28, 46, 6, 7], however, it has been observed that a change in the coupling strength is enough to create a chimera state by isolated desynchronization of all clusters except for one [32]. An example of one such representation of a chimera state can be found in Fig. 1(j).

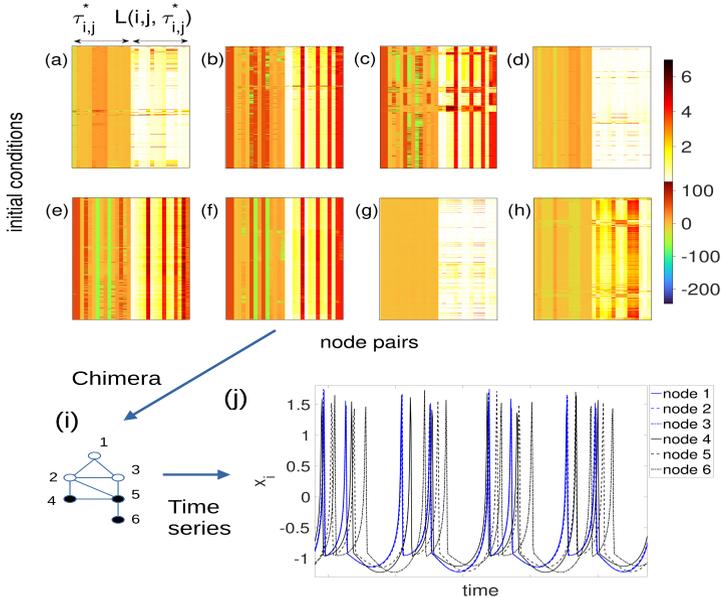

Supplementary Figure 1: All VPSs in Fig. 3 of the main text. (a)-(h) show all $e^l$ vectors (Eq. 12) in each of 8 clusters. The x-axis represents vertex pairs while the y-axis contains initial conditions forming a cluster. The color intensity represents entries of $e^l$. (i) The color-coded network shows the chimera state in (f). The vertices with a higher degree of synchronization, such that $L(i,j,\tau_{i,j}^*) \approx 0$, are colored white, while the remaining ones are shown in black. (j) shows $x_i$ versus time for this chimera state and is plotted as follows: first, $x_i(t)$ and $x_j(t-\tau_{i,j}^*)$ for nodes $i$ and $j$ having the smallest $L(i,j,\tau_{i,j}^*)$ value is plotted; thereafter, the phase-shift for any new node $k$ is decided such that $L(i,k,\tau_{i,k}^*)$ minimum over all the pre-existing nodes $i$.

The kind of synchrony present defines the details of the chimera state, whether the synchronous group is identical, generalized, phase synchronous, or approximate. In our synthetic model of HR oscillators on the DTI network as well as smaller models, we observe all of these. The stability of synchrony associated with a chimera state has some dependence on symmetry [32].

## 2.7 Nearly Synchronous States

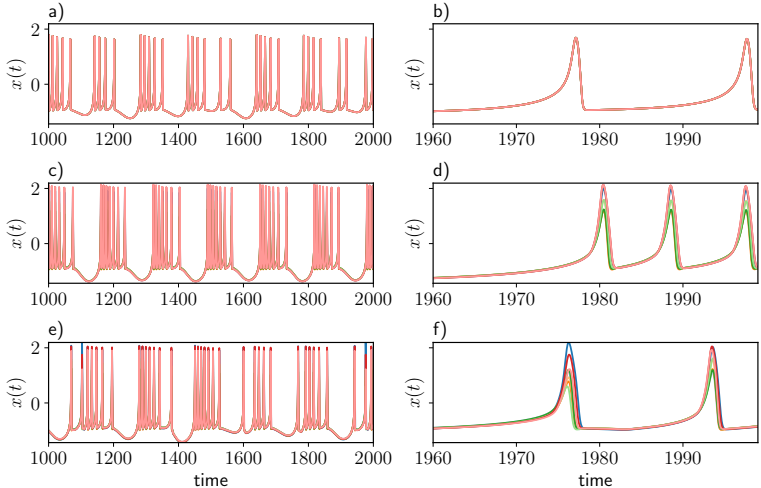

Supplementary Figure 2: Time series of the x-coordinate of the HR neuron model evolving on the DTI network without phase shift. (a) An x-coupled HR neuron model with no chemical coupling (i.e. $\alpha = 0$ and $\sigma = 1.7$), where we only observe identically synchronized oscillators. (b) Zoom in the panel of (a) for better visualization. (c) For chemical coupling $\alpha = 0.03$ (and $\sigma = 1.7$) oscillators still spike at the same time, but have different sized peaks creating a spread in the heights, though they still fully synchronize in phase. (d) Zoom in the panel of (c) for better visualization. (e) The HR model with the same coupling parameters as (c) but with the model parameter $a$ drawn from a normal distribution $a \sim \mathcal{N}(1, 0.1)$. (f) Zoom in the panel of (e) for better visualization. In this case, the oscillators peak in the same general time frame but these spikes have different peaks and different shapes. This represents a more generalized synchronization pattern.

Expanding upon the concept of approximate synchrony, Sec. 2.5, here for the HR model, we synonymously state as nearly synchronous. Consider Fig. 2 (a) and (b) in the model with only electrical coupling, Eq. 2 with dynamics $f$ given by Eq. 18 and coupling function $h_1$ from Eq. 19, yield identical synchronization as shown in Fig. 2 (a). However, in Fig. 2 (c) we have added a chemical coupling term changing the coupling function $h_2$ to the one from

Eq. 20, which is the model used by [18] (see Sec. 4). As can be seen, while the neurons still appear to be synchronized in phase, the peaks no longer match. In Fig. 2 (e) we have included both the electrical and chemical coupling terms, and additionally, we have introduced a slight parameter mismatch. The model parameter $a$, which was previously set to 1 for all oscillators is now drawn from a normal distribution $a \sim \mathcal{N}(1, 0.1)$. This is an example of nearly synchronous states of which is possible to examine the linear stability [44], even in the context of chimera [39].

## 3 Basin Structure and Fractals

### 3.1 Fractal Dimension

Examples of fractals are abundant in nature: they exhibit patterns across many scales that are self-similar, in some form. Real world examples, where the self-similarity is approximate and fine details do not continue indefinitely, include the coastlines of Britain (fractal dimension $\approx 1.21$) and Norway (fractal dimension $\approx 1.52$) [12]. Newton's fractals, Von-Koch curves, Cantor sets, Sierpinski's triangle, and the Koch Snowflake [12] are examples of mathematical fractals that continue indefinitely. Fractals can be understood as structures with a roughness that can be represented by a real number, called their dimension. A definitive property of fractals is the concept of non-integer dimensionality [12], as this fractional dimension motivates the very naming of sets as "fractals". Whereas fractals are all defined on the scaling of data density relative to an exponent describing the dimensionality, manifold dimensionality describes the integer number of coordinates necessary to uniquely parameterize a point. For instance, a straight line is 1 dimensional, a plane surface is 2 dimensional, and a cube is 3 dimensional. Fractal dimensionality is generally a description of the scaling behavior of a set with respect to decreasing window perspective; a non-integer fractal. It is a concept that applies to smooth as well as "rough" sets. From this perspective, the coastline of Norway is rougher than Britain, and correspondingly, its dimension is also farther from nearby integer numbers 1 and 2 [23, 5, 4, 40].

The fractal dimension is formally defined by the Hausdorff dimension, and various popular methods serve as useful estimates. For example, the similarity dimension, the correlation dimension, and the box dimension are each useful data-driven estimators of the fractal dimension with certain underlying assumptions on the data. For practical computational reasons, here we use the box dimension method to calculate the fractal dimension [35, 15, 37].

## 3.2 Box Dimension

The *box counting dimension* [15] is a popular estimator for fractal dimension. Squares (boxes or hyperboxes in higher dimensions) of side length $\epsilon$ are chosen and the number of 'boxes' ($N(\epsilon)$) of this size needed to cover the entire set are counted. The box dimension ($d$) is understood as the scaling exponent, in the limit as the area (or volume) of the box goes to zero with:

$$d_{\text{box}}(S_{BL}) = \lim_{\epsilon \to 0} \frac{\ln(N(\epsilon))}{\ln(1/\epsilon)}, \tag{5}$$

supposing the limit exists [43]. If the above limit does not exist, one may still find the upper and lower box dimensions of the set. Efficient algorithms for estimating the box dimension exist [37]. Fig. 3(a) shows an example of a boundary set obtained using Eq. 19, while Fig. 3(b) shows the plot of $\ln(\frac{1}{\epsilon})$ vs $\ln(N)$ for this boundary set. The slope of the line of best fit is interpreted as the limit given in Eq. 5, which in this example is equal to 1.27. As expected, the dimension of the boundary set is non-integer. Note that in practice due to finite data set size, as $\epsilon \to 0$, $N(\epsilon)$ saturates after a certain $\epsilon$ value is reached, which is due to the fact that $N(\epsilon)$ can not increase further once each point in the set is covered by a unique square.

## 3.3 Fractal Boundaries and Riddled Basins

Complexity can be found even in the structure of the basins of attraction in multi-stable systems. The feature of interest to this work is the possibility of the concept of a high level of "non-smoothness" of a fractal basin boundary when moving between sets of initial conditions attracted to one basin and those leading to a different one. In particular, we described a scenario that arises for the HR system and some parameters where a particularly highly intermingled fractal basin boundary called riddled basins occurs. We shall describe these scenarios.

In general, a boundary point of a set has the property that every neighborhood of the point, no matter how small, includes some points both in the set and not in the set. The boundary set of a given set is thus the union of all such boundary points. The boundary of an open disc in the plane, for example, is a circle, which clearly is a smooth set. However, in dynamical systems, it is not uncommon that boundary sets of basins of attraction can be fractals. The set of points forming the boundary between the basins of different attractors need not form a smooth curve, and instead, this boundary set may have a fractal dimension [26]. Our proposal herein is that the key feature of complex networked models of the brain with the observed behavior of rapid

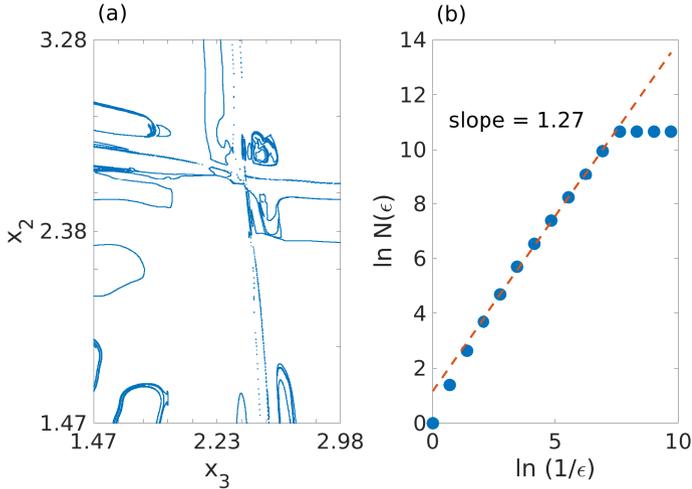

Supplementary Figure 3: (a) The boundary set obtained from Fig.3(b) in the main text. (b) The graph of $\ln(\frac{1}{\epsilon})$ vs $\ln N(\epsilon)$, calculating the dimension of the boundary set in (a). The fitted slope of 1.27 estimates the fractal dimension.

switching between various patterns is the presence of fractal basin boundaries [14, 26]. In particular, there is an apparently rich "intermingling" of these boundaries as the present phenomenon of what is called riddled basins [3, 29, 10, 30]; see Fig. 6. Riddled basins occur in many physical settings, notably in the basin structure, for example, of magnetic potential wells [47], but the fact that we point this out in neurophysiology as shown in Fig. 6, as a crucial mechanism behind the ability to quickly switch focus, the nimble brain theory, is a unique and new interpretation of this now classical concept from dynamical systems.

# 4 Computational Methods

## 4.1 Vector Pattern States

Here, we describe in more detail what we have introduced as Vector Pattern States (VPS) that allow comparisons between different initial conditions. The purpose of the VPS is to classify the various fully synchronous, chimera,

11

and incoherent states that may coexist in the basin. Thus, we associate which components of the complex network are in some form of synchrony, including the weaker approximate synchrony and phase-shifted synchrony as presented in Section 2. This is a crucial technological step for mapping the basins as described in the next section.

In general, synchrony is defined asymptotically (i.e. as $t \to \infty$), which is impractical. Instead for a finite time interval *time series*

$$\mathbf{x}_i(t), \text{ for } t \in [T_0, T_0 + T], i \in \{1, ..., N\}, \tag{6}$$

a similarity measure must be introduced to determine if nodes are 'close enough' to be considered synchronized. Let $I \subset \mathbb{R}$ and $\|\mathbf{u}\|_{L^2(I,M)}$ be the norm given by

$$\|\mathbf{u}\|^2_{L^2(I,M)} = \int_I \|\mathbf{u}(t)\|^2_2 dt, \tag{7}$$

where $M$ is a compact subset of $\mathbb{R}^d$ and $\|\cdot\|$ is the Euclidean norm $\|\mathbf{u}\|_2$ on $\mathbb{R}^d$. We abuse notation and also denote $\|\mathbf{u}\|_{L^2(I,M)}$ by $\|\mathbf{u}\|_2^2$. Also, we denote $u^\tau(t) = u(t - \tau)$ as the function $u$ evaluated at $t - \tau$. Other $L^p$ norms could be used to define the similarity measure, due to the equivalence of $L^p$ norms in finite-dimensional spaces. We will investigate the influence of different norms in future research.

An example of the similarity measure between the trajectories $\mathbf{x}_i(t)$ and $\mathbf{x}_j(t)$ of a pair $(i, j)$ could be an average over $I = [T_0, T_0 + T]$

$$L(i,j) = \frac{1}{T}\|\mathbf{x}_i - \mathbf{x}_j\|_2^2 = \frac{1}{T}\int_{T_0}^{T_0+T} \|\mathbf{x}_i(t) - \mathbf{x}_j(t)\|_2^2 dt. \tag{8}$$

This is fine in the CS case, however as we are interested in a more general synchronization case, the similarity measure should also allow for the possibility of a non-zero time shift $\tau$ between the time series. Thus, we consider the following similarity measure

$$L(i,j,\tau) = \frac{1}{T}\|\mathbf{x}_i - \mathbf{x}_j^\tau\|_2^2 = \frac{1}{T}\int_{T_0}^{T_0+T} \|\mathbf{x}_i(t) - \mathbf{x}_j(t - \tau)\|_2^2 dt. \tag{9}$$

The time shift that minimizes Eq.9, $\tau_{i,j}^*$, is of interest in this work since synchronization is viewed from a less restrictive criterion than CS, allowing for time shifts. For real-valued time series ($d = 1$), note that $L(i, j, \tau)$ is related to cross-correlation, in the sense that

$$(x_i \star x_j)_\tau = \int_{T_0}^{T_0+T} x_i(t)x_j(t-\tau)dt \tag{10}$$

is maximized at a particular $\tilde{\tau}$, which happens to be the same $\tau_{i,j}^*$ that minimizes $L(i,j,\tau)$. In other words,

$$\tau_{i,j}^* = \arg\max_{\tau} \ (x_i \star x_j)_\tau = \arg\min_{\tau} \ L(i,j,\tau). \qquad (11)$$

Thus, it is possible to find $\tau_{i,j}^*$ using the well-studied method of cross-correlation. In practice, not only is the time interval finite, but the sampling rate is finite as well, allowing the integrals to be replaced by sums and the cross-correlation $R_{x_i,x_j}(\tau)$ can be found using the xcorr function in MATLAB. Since we define $\tau_{i,j}^* = \arg\max_{\tau} R_{x_i,x_j}(\tau)$, when $T - |\tau_{i,j}^*| > 0$, we calculate $L(i,j,\tau_{i,j}^*)$ using the formula

$$L(i,j,\tau_{i,j}^*) = \frac{1}{T - |\tau_{i,j}^*|} \sum_{t=T_0+|\tau_{i,j}^*|}^{T_0+T} \|\mathbf{x}_i(t) - \mathbf{x}_j(t - \tau_{i,j}^*)\|_2^2.$$

In our simulations, the time interval and sampling rate are chosen so that $R_{x_i,x_j}(\tau)$ gives a good approximation of $\tau_{i,j}^*$. For instance, for the HR system, we select the fast variable ($x$-coordinate), so the time series has a sufficiently large number of cycles for the chosen time interval and step size.

For the case of phase synchronization (where only phase difference is relevant), the $\tau_{i,j}^*$ for every $i \neq j$ encodes all information about different final states, that is, they distinguish between chimera, complete synchronous, and completely incoherent states. In the generalized version of synchronization discussed here, not only are the $\tau_{i,j}^*$ important but so is the difference between the 'profile' of the time series, so $L(i,j,\tau_{i,j}^*)$ should be taken into consideration as well. Therefore the state of the pair $(i,j)$ can be further characterized by the ordered pair $(\tau_{i,j}^*, L(i,j,\tau_{i,j}^*))$.

We are now ready to define the Vector Pattern State (VPS) for all $(i \neq j)$

$$e^l = (\tau_{1,2}^{*,l}, \tau_{1,3}^{*,l}, \ldots, \tau_{N-1,N}^{*,l}, \beta L^l(1,2,\tau_{1,2}^{*,l}), \beta L^l(1,3,\tau_{1,3}^{*,l}), \ldots, \beta L^l(N-1,N,\tau_{N-1,N}^{*,l}), \qquad (12)$$

where $\beta > 0$ allows for scaling the relative importance of the phase shift vs. the similarity measure between the time series. Each initial condition will lead to its own VPS. Though such states may not be unique, for instance, an initial condition that leads to CS will have the VPS: $e^l = (0,0,0,\ldots,0)$ (in the limit as $t \to \infty$), and thus each initial condition that leads to CS will have the same VPS. In the above mentioned case of phase synchronization, the VPS will be given by $e^l = (\tau_{1,2}^*, \tau_{1,3}^*, \ldots, \tau_{N-1,N}^*, 0, 0, \ldots, 0)$ (in the limit as $t \to \infty$).

## 4.2 Mapping the Basin Structure

In the last section, we described our VPS so as to compare components of network coupled dynamics. With this, we can map the entire space of initial conditions to decide on a basin structure. The VPS allows for the characterization of states that may have both a phase shift as well as differences in the 'shape' of the time series. As discussed above, in the CS scenario, all components of the VPS will be $\approx 0$. Similarly in the phase synchronization, all components corresponding to $L(i, j, \tau_{i,j}^*)$ will be 0, while the $\tau_{i,j}^*$ themselves become bounded from above. However, the utility of the VPS is best expressed in chimera states and generalized synchronization. When some of the vertices become asynchronous and others follow a CS state, the VPS will have some components nearly 0 and others large.

We examine the CS case in more detail, noting it is analogous to the more generalized synchronization. Finding the basin structure is straightforward if the pairwise distances (under whatever measure is chosen) between all oscillators are small enough. The system may be considered in CS and the initial condition leading to this state may be marked as one color. Any initial condition not leading to CS will then be assigned a different color. However, this becomes much more challenging in chimera states, since not all pairs will reach CS. It may seem that classifying all states that have at least one synchronized pair as well as at least one desynchronized pair in the same basin will be sufficient, but chimera states range from being nearly incoherent to being nearly coherent, and classifying all as the same may hide very rich structure.

To resolve this, each VPS ($e^l$) is stacked into a matrix,

$$E = \begin{bmatrix} e^1 \\ e^2 \\ \vdots \\ e^S \end{bmatrix}, \qquad (13)$$

where $S$ is the number of VPSs or sample size. This matrix encodes all final states evolving from each initial condition, which are indexed by $l$. $E$ is then clustered, in this case using $k$-means (see Section 4.3), with each cluster assigned a color in the basin plot. This leaves the choice of $k$, which was done using a variation on the elbow method. An example of how this is done is shown in Fig. 5.

In the elbow method, after a cluster has been chosen from $k$-means, the value of $W(C)$ from Eq. 16 becomes the dependent variable, with $k$ as the independent variable. When examining this relationship across multiple $k$'s

14

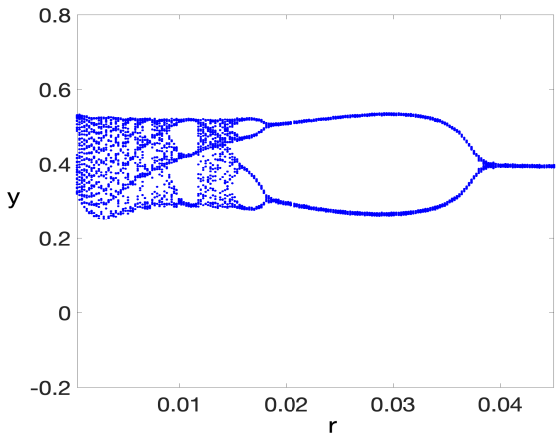

Supplementary Figure 4: The bifurcation diagram for an isolated HR oscillator. The parameters are $a = 1, b = 3, c = 1, d = 5, x_R = -0.5(1 + \sqrt{5})$, $s = 4$, and $I = 3.27$ [11].

generally there will be a discontinuity that is clearly visible, which will be the optimal value, though there are situations where this can fail [16]. In a slight variation, to identify the elbow, a log-log plot is used instead, and the last $k$ before the majority of values fall on the linear trend is selected as the 'elbow' as shown in Fig. 5 (b).

The effects varying $k$ has on the basin structure are shown in Fig. 6, the initial conditions chosen from the random plane shown here were drawn from a $1,500 \times 1,500$ grid, sampled over the unit square ($[0,1]^2$). The initial conditions for all other oscillators were initially chosen uniformly at random from the uniform distribution $\mathcal{U}(-1, 1)$ and were then held constant with the exception of the initial conditions drawn from the random plane (which were sampled as discussed above). The complexity of the basin can be seen even for small values of $k$, the optimal $k$ balances between too little and too much detail in the basin structure. In Fig. 6, $k = 8$ was found to be optimal using the version of the elbow method discussed above. Note that the basin used in this example differs from the one shown in Fig. 5, hence the different values of $k$.

15

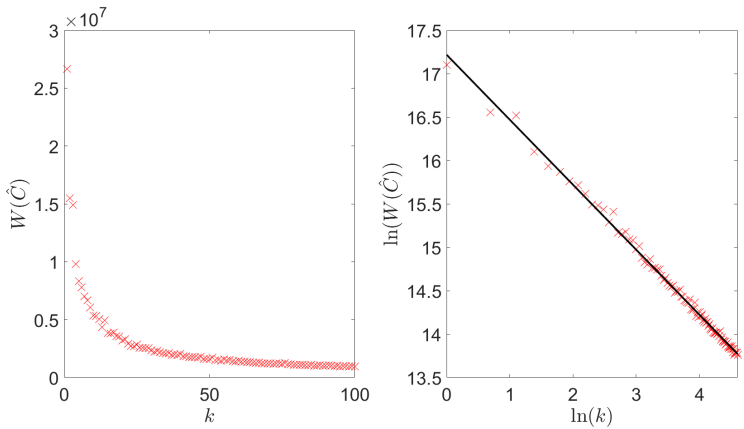

Supplementary Figure 5: The elbow method. (a) The traditional version of the elbow method relies on examining the plot of $W(\hat{C})$ and finding a point discontinuity, which is then chosen as $k$. Here it is challenging to find a clear elbow. (b) Instead, $\ln(k)$ vs. $\ln(W(\hat{C}))$ is plotted, with a linear fit determined. The best $k$ value is chosen to be the last value before which the majority of the points fall on the line of best fit, in this example $k = 5$.

### 4.3 Clustering and k-Means

Often insights into a problem can be drawn by clustering data together to infer meaningful relationships. Numerous methods exist to perform clustering, however, k-means is a popular, robust, and relatively computationally inexpensive method for clustering a 'cloud' of data points together. With data $x \in \mathbb{R}^{S \times d}$, where $S$ is the sample size and $d$ is the dimension, and making the assumption that data is all of the quantitative types so that distance has meaning, then k-means is performed by minimizing the within cluster mean squared Euclidean distance between points. Choosing $1 \leq k \leq S$, let $C$ be a clustering $C = \{C_1, C_2, ..., C_k\}$, with $C_k \subset \{1, 2, ..., S\}$. Then we define the dissimilarity measure $d$ to be the squared Euclidean distance [16]:

$$d(x, y) = ||x - y||_2^2, \tag{14}$$

16

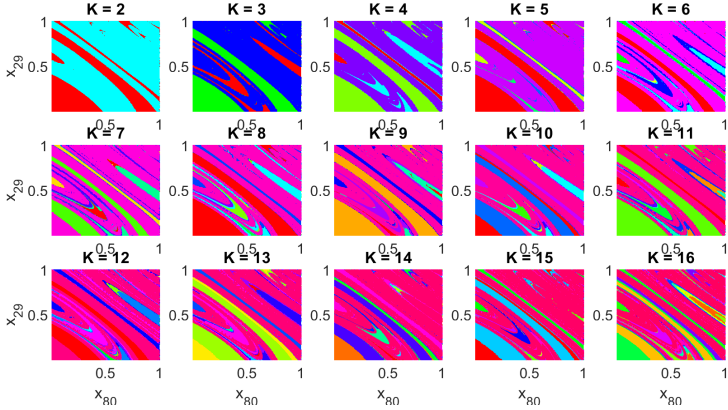

Supplementary Figure 6: Effect of varying $k$ on the basin plots. Even for small $k$, a complex basin structure is observed. The choice of $k$ can be seen as a trade off between too much and too little detail in this context.

and the within cluster mean

$$\bar{x}_k = \frac{1}{N_k} \sum_{i=1}^{S} I(i \in C_k) \cdot x_i, \tag{15}$$

where $N_k = \sum_{i=1}^{S} I(i \in C_k)$ and $I$ is the indicator function. The within cluster squared error $(W)$ is defined as,

$$W(C) = \sum_{j=1}^{k} \sum_{i=1}^{S} I(i \in C_j) d(x_i, \bar{x}_j). \tag{16}$$

Note that $d(x_i, \bar{x}_j)$ is the distance between a data point and the centroid of the cluster it is in. The k-means problem is to minimize this within cluster squared error, that is to find the cluster $\hat{C}$ with

$$W_k(\hat{C}) = \min_{C} W(C). \tag{17}$$

As noted above, there are numerous clustering techniques, and one can even generalize k-means by designing the distance function $d(x, y)$ to be a different dissimilarity measure than the more common squared Euclidean distance, such a clustering method is called k-medoids [16].

# 5 Networks of coupled oscillators

Here, we show that the construction of basins of attractions via VPS can also be applied to networks of coupled oscillators. First, we give details about the network dynamics used in the main text. Subsequently, we illustrate our approach to other network dynamics: Kuramoto oscillators and Hénon maps. In all cases our approach is able to capture the basin structure, illustrating its applicability.

## 5.1 Coupled HR neurons

Hindmarsh and Rose (HR) developed a model of neuronal firing in the brain [17]. We have used this model as a complex network of coupled oscillators as a semi-synthetic system using a true DTI network. Each HR neuron follows the following dynamics:

$$f(\mathbf{x}) = \begin{bmatrix} y - ax^3 + bx^2 - z + I \\ c - dx^2 - y \\ r[s(x - x_R) - z], \end{bmatrix} \quad (18)$$

where $\mathbf{x} = [x, y, z]^T$. The parameters $a$, $b$, $c$ and $d$ are all polynomial parameters, generally they are chosen to be $a = 1, b = 3, c = 1, d = 5$ following [17] and in this work we do the same. Varying $b$ allows for switching between bursting and spiking behavior, $I$ mimics the membrane input potential, $r$ controls how quickly the slow variable ($z$) varies, $s$ governs adaptation, and $x_R$ represents the resting potential. Additionally, two different coupling functions were used, the first one,

$$h_1(\mathbf{x}_i, \mathbf{x}_j) = \begin{bmatrix} x_j - x_i \\ y_j - y_i \\ z_j - z_i \end{bmatrix}, \quad (19)$$

and the other one with

$$h_2(\mathbf{x}_i, \mathbf{x}_j) = \begin{bmatrix} 0 \\ y_j - y_i \\ 0 \end{bmatrix} - \alpha(x_i - V_{syn}) \begin{bmatrix} [1 + e^{-\lambda(x_j - \theta_{syn})}]^{-1} \\ 0 \\ 0 \end{bmatrix} \quad (20)$$

For the standard coupling scheme with $h_1$ given in Eq. 19, the values of $a, b, c, d, x_R$, and $s$ are selected based on earlier research works, see [11, 17]. In this work, we set these parameters to $x_R = -0.5(1+\sqrt{5})$, $s = 4$, and $I = 3.27$. This will be referred to as the standard HR model. The value of $x_R$ is chosen from the condition $\dot{x}_i = \dot{y}_i = 0$ while $z_i$ (the slow variable) is neglected. With

this choice of parameters, and depending on $r$, an isolated HR neuron can exhibit a periodic time evolution. To decide if the time evolution is periodic or not, we follow the approach of Poincare [43]. Fig. 4 highlights how the behavior of the system changes when the parameter $r$ is varied. For a given set of parameters, the number of cuts a moving HR neuron makes on the Poincaré plane is counted. The Poincaré plane is any plane transverse to a time-continuous flow, which in this case is the $y - z$ plane located at $x = 0$. Fixing $I$ and decreasing $r$ towards 0, the $y$ coordinates of all intersections on the Poincare plane are shown, such that $x(t) < 0$ and $x(t + dt) > 0$, where $dt$ denotes the integration time step. The resulting plot (Fig. 4) is called a bifurcation diagram and the type of bifurcation is a period-doubling bifurcation of limit cycles. A finite number of cuts (let's say $N_c$) indicate a periodic trajectory with period $N_c$, while the limit $N_c \to \infty$ shows a chaotic evolution.

An alternative model, which includes chemical coupling in addition to the electrical coupling, that is with $h_2$ from Eq. 20, the parameters are: $a = 1, b = 3, c = 1, d = 5, s = 4, r = 0.005, x_R = -1.6, I = 3.25, \sigma = 0.5, \alpha = 0.03, V_{syn} = 2, \theta_{syn} = -0.25$ and $\lambda = 10$, which is known to contain chimeras [18]. Here the chemical and electrical adjacency matrices are assumed to be the same, unlike in [18]. This will be referred to as the chemical coupling model. An example of the dynamics is shown in Fig. 2 (c) and (d), and a version of this model where the parameter $a$ for each oscillator is drawn from a normal distribution is shown in Fig. 2 (e) and (f).

## 5.2 Two populations of globally coupled oscillators

Here, we consider the case of two small size populations that are globally coupled to each other [1]. We consider two coupled populations as illustrated in Figure 7 a)

$$\dot{\theta}_i^k = \omega + \sum_{k'=1}^{2} \frac{\sigma_{kk'}}{N_{k'}} \sum_{j=1}^{N_{k'}} \sin(\theta_j^{k'} - \theta_i^k - \alpha), \quad i = 1, \ldots, N_k, \qquad (21)$$

where $\theta_i^k$ is the phase of the $i$-th oscillator in population $k \in \{1, 2\}$ and $\omega$ is the oscillator frequency. Based on earlier works [1, 24, 31], we suppose that $\sigma_{11} = \sigma_{22} = \mu > 0$, and $\sigma_{12} = \sigma_{21} = \nu > 0$, with $\mu > \nu$. Also, by rescaling time, we may set $\mu + \nu = 1$ and introduce useful parameters $A = \mu - \nu$ and $\gamma = \pi/2 - \alpha$ that determine the existence of chimeras states in the system. For identical oscillators, we can change coordinates to the rotating frame $\omega$, $\theta_i^k \mapsto \theta_i^k - \omega t$, so the natural frequency can be fixed at $\omega = 0$.

19

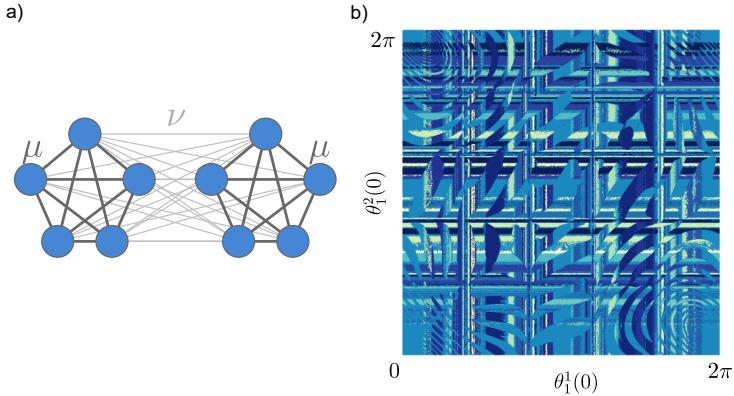

Supplementary Figure 7: Basin structure of two globally coupled populations. (a) Two globally coupled populations of five nodes. The (intra-)coupling strength inside a population and (inter-)coupling strength between populations are given by $\mu$ and $\nu$, respectively, with $\mu > \nu$. (b) Two-dimensional section of the state space showing basins of the 18 (clustered) distinct states identified by the VPS of the two globally coupled populations. There is a symmetry among the basins with respect to reflections across the diagonals, which originates from the reflection symmetry of the network. To construct the VPS, we use $\beta = 1$ in Equation (12) and a grid of $1,257 \times 1,257$ uniformly sampled initial conditions.

The system (21) has full permutation symmetry, allowing a low-dimensional description in terms of macroscopic variables in the thermodynamic limit ($N_k \to \infty$) [1]. In particular, the authors in [24] use the resulting reduced system to describe the basins of attraction of chimera states. The basin structure corresponds to the infinite size system, but it is expected to remain roughly the same for large population sizes [31]. We tested our VPS construction (results not shown), reproducing a basin structure that contains a few distinct states, as predicted by the infinite size system [24]. Here, we consider the interesting case when the population size is small, and the dynamics of the full system deviate from the evolution of the reduced system [31]. In fact, in this scenario, the evolution still has a low-dimensional description but the macroscopic variables do not evolve in a closed-form equation, so the reduced system becomes intractable analytically, see [25] for details.

We consider the case of $N = 5$, $A = 0.1$ and $\gamma = 0.025$, where it is known that stable chimeras exist [1, 24, 31]. For initial conditions we consider equally spaced phases over the circle, i.e., for a fixed $\vartheta \in [0, 2\pi]$

$$\theta_i^1(0) = \frac{2\pi i}{N}, \quad i = 1, \ldots, N_1,$$

$$\theta_i^2(0) = \frac{2\pi i}{N} + \vartheta, \quad i = 1, \ldots, N_2,$$

so the system starts close to an incoherence state.

We consider $\vartheta = 0.1$, transient time of $2,000$ time steps, and integrate over 1000 more time steps to obtain the final state, sampling every 0.02 time step. The basin plot is performed in the plane $(\theta_1^1(0), \theta_1^2(0))$ for which the initial conditions are varied over the interval $[0, 2\pi]$, see Figure 7 (b). We see that the basin has a coexistence of different chimera states and (complete) phase synchronization states, as described in [31], and our method is capable of identifying these cases.

## 5.3 Coupled oscillators in networks

Here, we apply the VPS construction to networks that are not globally coupled, but rather have some missing links. In this scenario, instead of the network allowing all permutation symmetries, only a few permutations are possible. We consider the following equations of motion for the oscillators

$$\dot{\theta}_i = \sigma \sum_{j=1}^{N} [A]_{i,j} \sin(\theta_j - \theta_i - \alpha), \quad i = 1, \ldots, N, \qquad (22)$$

where $\sigma = 0.01$ is the overall coupling strength. The adjacency matrix $A$ represents a network that does not have full permutation symmetry. To generate this network we initiate the globally coupled network, as depicted in Figure 7 (a), and remove uniformly at random one edge from the graph. The resulting graph has two new clusters that are depicted in different colors, which is calculated via Sage [41]. See Figure 8 for an illustration of the (random) network used in our numerical results.

We consider a time interval of $5,000$ time steps and sample the trajectory at every 0.02 time step. Also, the parameters are $\sigma = 0.01$ and as before $\alpha = \pi/2 - \gamma$ with $\gamma = 0.025$. Choosing the same two-dimensional section as previously, we construct the basin structure via VPS, as illustrated in Figure 8 (b), using five clusters that are determined by the elbow method discussed earlier. Here, we consider a transient time of $4,000$ time steps and calculate

21

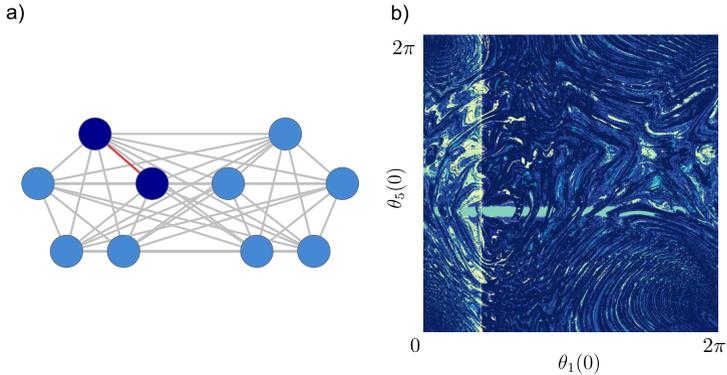

Supplementary Figure 8: (a) Network with one single link removed, breaking the full permutation symmetry of the globally coupled network. Nodes with different colors correspond to the new clusters formed from the permutation symmetry group. The red edge corresponds to a randomly selected edge that was removed from the globally coupled network. (b) Two-dimensional section of the state space showing basins of 12 (clustered) distinct states. To construct the VPS, we use $\beta = 1$ in Equation (12) and a grid with $1,257 \times 1,257$ uniformly sampled initial conditions.

the similarity measure and cross-correlation, over the remaining $1,000$ time steps.

We also test the performance of the method for another network structure, where two links are removed, see Figure 9 a). We consider a time interval of $15,000$ time steps, discarding the first $14,000$ time steps as transient, and sample the trajectory at every $0.05$ time step. The parameters are $\alpha = \pi/2 - \gamma$ with $\gamma = 0.075$ and $\sigma = 0.01$. Figure 9 b) displays the basin structure identified by our approach, where a zoom region has been enlarged, confirming that the basin is intermingled.

In both cases, the basins show that the systems have a high sensitivity to initial conditions, so our method is able to characterize such basins that have intermingled structures. It remains an interesting line of research to describe the mechanism behind the chaotic dynamics since it differs from the globally coupled scenario where a low-dimensional description is available [8].

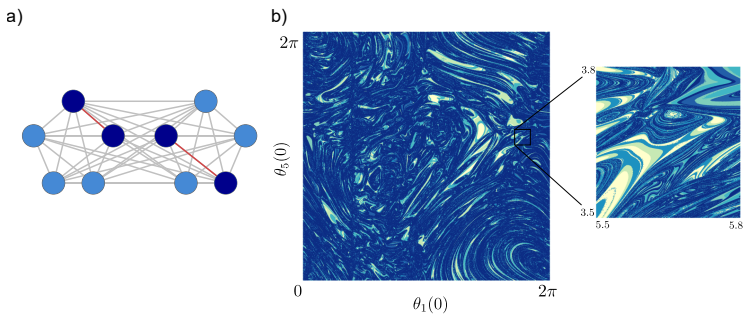

Supplementary Figure 9: (a) Network with two links removed, breaking the full permutation symmetry of the globally coupled network. Nodes with different colors correspond to the new clusters formed from the permutation symmetry group. The red edges correspond to randomly selected edges that were removed from the globally coupled network. (b) Two-dimensional section of the state space showing basins of 7 (clustered) distinct states. To construct the VPS, we use $\beta = 1$ in Equation (12) and a grid with $1,257 \times 1,257$ and $624 \times 624$ for the left and right panel, respectively, uniformly sampled initial conditions. The parameters are $\gamma = 0.075$ and $\sigma = 0.01$.

## 5.4 VPS for phase oscillators

Different from the HR system, when we consider coupled oscillators, the time series is given by phases $\{\theta_i(t)\}_{t\geq 0}$ in the circle. To construct the VPS, we perform an embedding into $\mathbb{R}$, so we can use the same similarity measure introduced in Section 4.1. More precisely, for each oscillator time series, we consider

$$x_i(t) = \cos(\theta_i(t)), \quad i = 1, \ldots, N.$$

## 5.5 Coupled Hénon maps

To allow the system to settle on to the attractor we use $20,000$ time steps. The final $1,000$ time steps are utilized for the construction of the VPS. Following several of the above examples, we use only the x-component for the VPS. The basin shown contains $1,000,000$ initial conditions sampled in a grid of $[-3,3]^2$ for the $x$ and $y$ component of a randomly chosen node. The initial condition for all other nodes is set to $(0,0)^T$. Though we are using a different network structure, we find a very similar basin structure to that

of [36]. This basin appears to be riddled, which indicates how general this phenomenon is. Additionally, it is clear that non-trivial symmetry is not necessary for riddled basins.

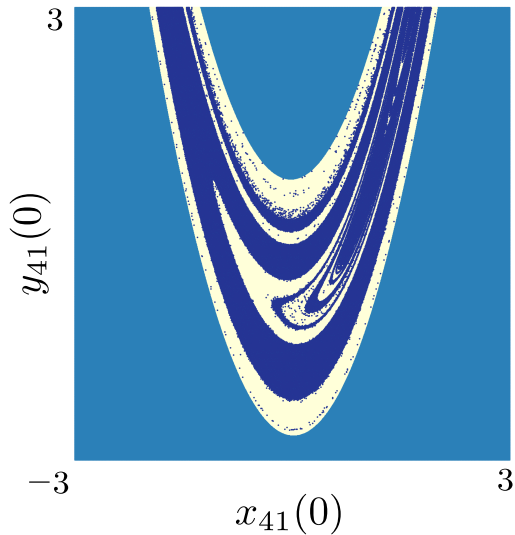

Supplementary Figure 10: Basin structure in the Hénon system for a randomly chosen x,y plane. The network is the same DTI network shown in Fig.1 in the main text. This structure appears to have a riddled basin.

 
\end{document}